\documentclass[10pt, twocolumn, twoside]{IEEEtran}

\usepackage{pifont}
\usepackage{cite}
\usepackage[dvips]{epsfig}
\usepackage{graphicx}
\usepackage{calligra}
\usepackage[T1]{fontenc}
\usepackage{bbold}
\usepackage{mathrsfs}

\usepackage{subfigure}
\usepackage{color}
\usepackage{amsmath}
\usepackage{cases}

\usepackage{amssymb}
\usepackage[justification=centering]{caption}
\usepackage{framed}

\usepackage{subfigure}
\usepackage{color}
\usepackage{amsmath}
\usepackage{bm}
\usepackage{amssymb}
\usepackage{amsbsy}
\usepackage{bbm}
\usepackage{dsfont}
\usepackage[normalem]{ulem}
\usepackage{xcolor,cancel}
\usepackage{enumerate}
\usepackage{paralist}

\newtheorem{theorem}{Theorem}
\newtheorem{proposition}{Proposition}
\newtheorem{lemma}{Lemma}

\newtheorem{assumption}{Assumption}
\newtheorem{definition}{Definition}

\ifCLASSINFOpdf
\else
\fi
\begin{document}
%
\title{A Fundamental Limitation on Maximum Parameter Dimension for Accurate Estimation with Quantized Data 
}
%
%
%

\author{Jiangfan Zhang,
	Rick S. Blum,~\IEEEmembership{Fellow,~IEEE}, Lance Kaplan,~\IEEEmembership{Fellow,~IEEE}, and Xuanxuan Lu

\thanks{This work was supported by the U. S. Army Research Laboratory and the U. S. Army Research Office and was accomplished under Agreement Numbers W911NF-14-1-0245 and W911NF-14-1-0261. The views and conclusions contained in this document are those of the authors and should not be interpreted as representing the official policies, either expressed or implied, of the Army Research Laboratory, Army Research Office, or the U.S. Government. The U.S. Government is authorized to reproduce and distribute reprints for Government purposes notwithstanding any copyright notation here on.}



}


\maketitle

\begin{abstract}

It is revealed that there is a link between the quantization approach employed and the dimension of the vector parameter which can be accurately estimated by a quantized estimation system. A critical quantity called inestimable dimension for quantized data (IDQD) is introduced, 
which doesn't depend on the quantization regions and the statistical models of the observations but instead depends only on the number of sensors and on the precision of the vector quantizers employed by the system. It is shown that the IDQD describes a quantization induced fundamental limitation on the estimation capabilities of the system. To be specific, if the dimension of the desired vector parameter is larger than the IDQD of the quantized estimation system, then the Fisher information matrix for estimating the desired vector parameter is singular, and moreover, there exist infinitely many nonidentifiable vector parameter points in the vector parameter space. Furthermore, it is shown that under some common assumptions on the statistical models of the observations and the quantization system, a smaller IDQD can be obtained, which can specify an even more limiting quantization induced fundamental limitation on the estimation capabilities of the system. 
\end{abstract}


\begin{IEEEkeywords}
	
	Distributed sensor parameter estimation, inestimable dimension for quantized data, singular Fisher information matrix, identifiability, quantization.
	
\end{IEEEkeywords}

%

\section{Introduction}

Bolstered by recent technological advances in coding, digital wireless communications technology and digital devices, the employment of quantized data has become increasingly popular in many applications, such as sensor networking, the internet of things, data-transmission systems and  data-storage systems. Inspired by this tendency, parameter estimation utilizing quantized data has seen considerable interest in recent years, see \cite{niu2006target, chen2010nonparametric, kar2012optimal, shen2014robust, ozdemir2009channel, fang2009hyperplane, venkitasubramaniam2007quantization,    zhang2015Asymptotically} and references therein. 

\begin{figure}[htb]
	\centerline{
		\includegraphics[width=0.46\textwidth]{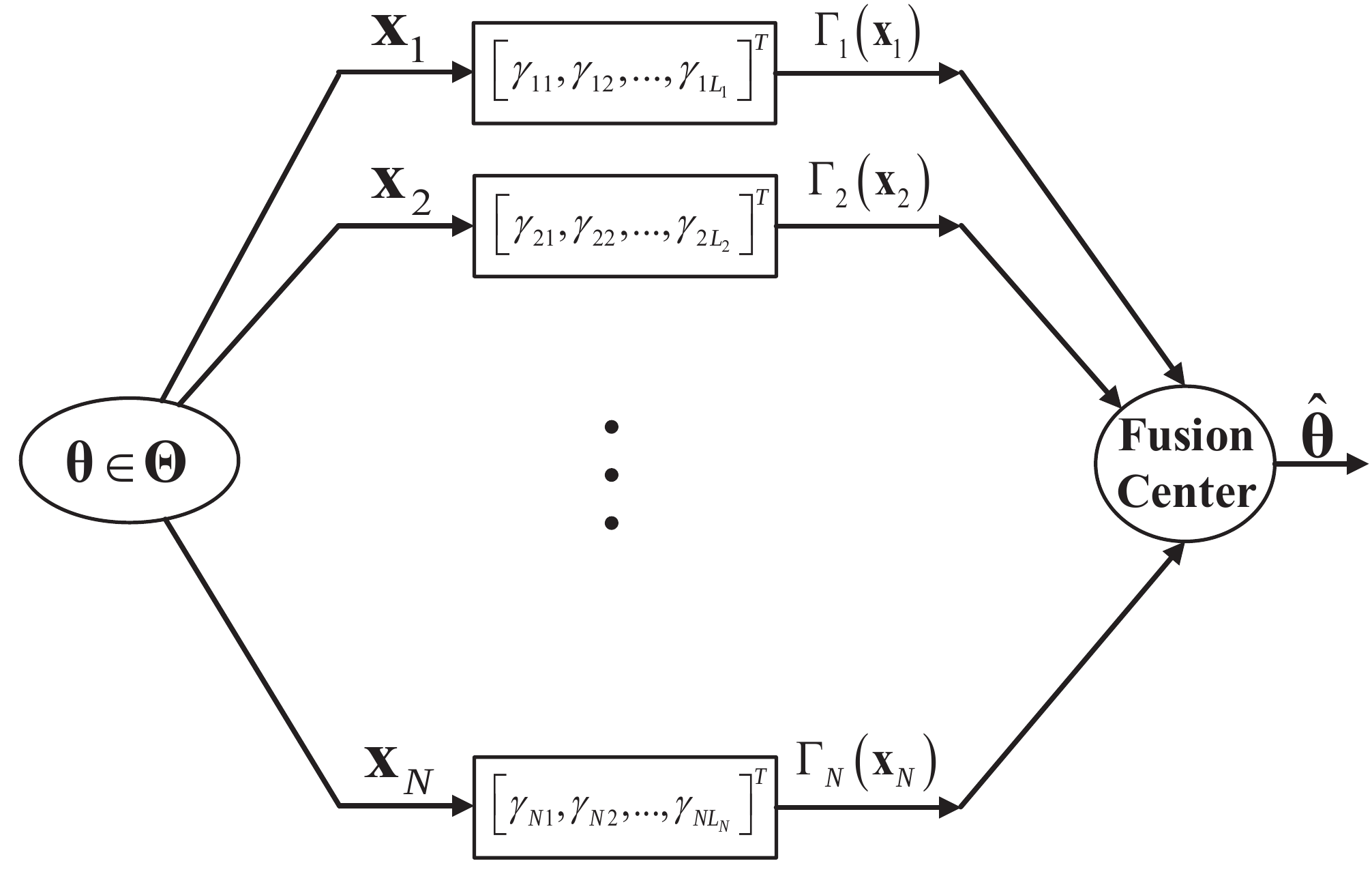}
	}
	\caption{Parameter estimation system with quantized data.}
	\label{Fig_estimation_system}
\end{figure}

\begin{figure}[htb]
	\centerline{
		\includegraphics[width=0.46\textwidth]{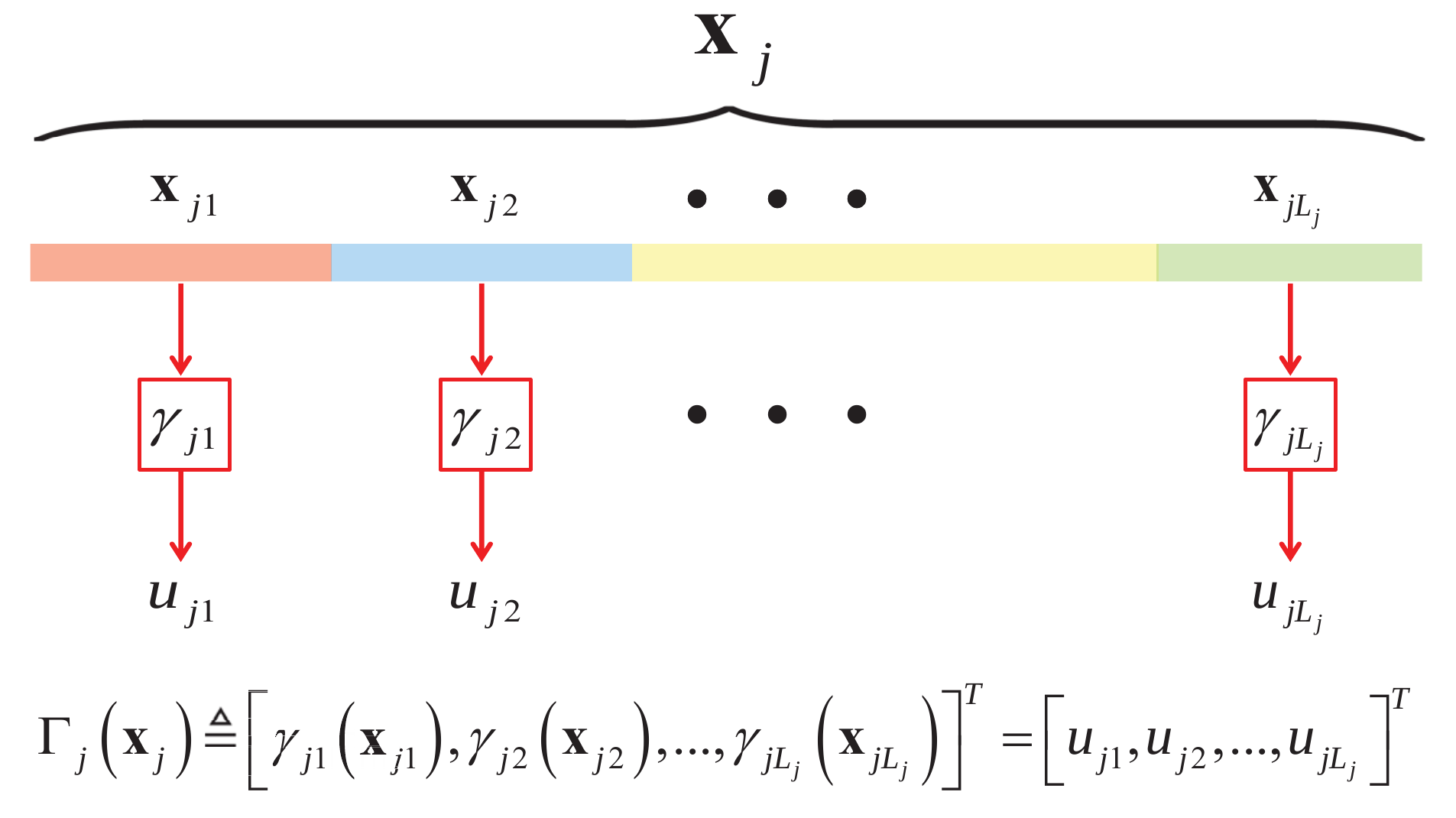}
	}
	\caption{The sequence of vector quantizers employed for ${\bf x}_j$.}
	\label{Fig_quantizer}
\end{figure}

A parameter estimation system which employs quantized data is depicted by Fig. \ref{Fig_estimation_system}. The  distribution of the observations $[{\bf x}_1^T, {\bf x}_2^T, \cdots, {\bf x}_N^T]^T$ depends on an underlying vector parameter $\boldsymbol{\theta} \in {\boldsymbol{\Theta}}$. As shown in Fig. \ref{Fig_estimation_system}, for each $j$, a sequence of $L_j$ vector quantizers denoted as ${\Gamma _j} \buildrel \Delta \over = {[ {{\gamma _{j1}},{\gamma _{j2}},...,{\gamma _{j{L_j}}}} ]^T}$ is employed to convert the observation vector ${\bf x}_j$ to digital data ${\Gamma}_j({\bf x}_j)$, which is transmitted, without error,  to the fusion center (FC). To be specific, as illustrated by Fig. \ref{Fig_quantizer}, ${\bf x}_j$ is partitioned into a sequence of $L_j$ disjoint observation subvectors $\{{\bf x}_{jl}\}_{l=1}^{L_j}$ first, and then for each $l$, ${\bf x}_{jl}$ is quantized to $u_{jl}$ by the $l$-th vector quantizer $\gamma_{jl}$ in the sequence ${\Gamma _j}$. The output of the sequence of vector quantizers ${\Gamma}_j({\bf x}_j)$ is the vector $[u_{j1},u_{j2},...,u_{jL_j}]^T$ which gathers the quantized data from all vector quantizers $\{\gamma_{jl}\}$.
After collecting the quantized data from all sequences of vector quantizers, the FC makes use of $\{{\Gamma}_j({\bf x}_j)\}$ to estimate the value of the desired vector parameter $\boldsymbol{\theta}$. 

In general, the output of a quantizer is a scalar, however, the output of ${\Gamma _j}$ is a vector. To distinguish ${\Gamma _j}$ from the commonly defined quantizers, we refer to the sequence of vector quantizers, ${\Gamma _j}$, as a superquantizer. It is worth mentioning that $L_j$ can be any positive integer. The scenarios where $L_j>1$ are widely considered in recent literature, see \cite{shen2014robust} for instance. For the scenario where $L_j=1$,  the superquantizer ${\Gamma _j}$ is equivalent to a vector quantizer.

Our recent investigations into attacks on parameter estimation systems provide insight into some very effective attacks on systems utilizing quantized data \cite{zhang2015Asymptotically}. Building on these ideas, this paper attempts to fully uncover the fundamental limitations on the estimation capabilities of the unattacked quantized estimation system shown in Fig. \ref{Fig_estimation_system}. In order to assess the estimation capabilities of the quantized estimation system shown in Fig. \ref{Fig_estimation_system}, two criteria are often adopted \cite{hochwald1997identifiability}. The first criterion is the information-regularity condition, which is defined as
\begin{definition}[Information-Regularity Condition]  \label{Definition_Information_Regularity}
	The Fisher information matrix (FIM) for estimating the desired parameter is nonsingular.
\end{definition}  

The information-regularity condition guarantees the existence of the Cramer-Rao bound (CRB) for the desired parameter. Further, under mild additional conditions it guarantees the estimation performance of an appropriate estimator can always be improved by an increase in the number of observations provided a suitably large set of observations is employed. 
Moreover, it can be shown that in most of cases, if the FIM is singular, there is no unbiased estimator for the desired parameter with finite variance \cite{stoica2001parameter}.
To this end, the information-regularity condition, which ensures the nonsingularity of the FIM for the desired parameter $\boldsymbol{\theta}$,
is crucial in parameter estimation problems.  

The second criterion to assess the estimation capabilities of the quantized estimation system is the identifiability condition, which is defined as
\begin{definition}[Identifiability Condition]\label{Definition_Identifiability}
	There exists no parameter value in the parameter space such that the conditional distribution of the data conditioned on the parameter is identical to that for some other parameter value in the parameter space. 
\end{definition}
  
The identifiability condition is sufficient to guarantee almost sure convergence of a class of estimators which includes the maximum likelihood estimator as the sample size approaches infinity, given some reasonable conditions \cite{wald1949note, hochwald1997identifiability}. 
Intuitively, if the identifiability condition fails for some parameter value, then there exists another parameter value which is just as likely as the true value based on the observations,
and hence the desired parameter cannot be estimated consistently \cite{lehmann1999elements, bekker2001identification}. Therefore, it is necessary that every parameter point in the parameter space ${\boldsymbol{\Theta}}$ satisfies the identifiability condition for a meaningful estimation problem. In such cases, we say the vector parameter space is identifiable.

In general, the FIM nonsingularity and the identifiability of the vector parameter space are both determined by the statistical models of the observations, the value of the desired parameter, and the quantizer designs employed by the system which complicates the analysis of the estimation capabilities of a quantized estimation system.  On the other hand, this paper provides a simple expression called the inestimable dimension for quantized data (IDQD) which describes a vector parameter dimension beyond which the FIM nonsingularity and the identifiability of the vector parameter space are both guaranteed to fail. These powerful results can be employed in preliminary design in many applications of quantized (digital) data and seem especially important for big data problems which are attracting significant attention lately.

\subsection{Summary of Results}
\begin{asparaenum}
\item For the general parameter estimation with quantized data system shown in Fig. \ref{Fig_estimation_system}, the impact of quantization with regard to the information-regularity condition is first studied. 
By exploring the structure of the FIM for estimating the desired vector parameter, 
it is shown that if the dimension of the desired vector parameter is larger than the IDQD, the FIM for estimating the desired vector parameter cannot be nonsingular for any statistical models of the observations, any value of the desired vector parameter, and any quantization regions. 
Hence, the IDQD specifies a quantization induced fundamental limitation on the estimation capabilities of the quantized estimation system with regard to the information-regularity condition, which limits the number of parameters which can be estimated by the quantized estimation system while maintaining a nonsingular FIM.  

\item Next, we investigate the impact of quantization on the identifiability condition. It is shown that for any statistical models of the observations and any quantization regions employed by the system,  if the dimension of the desired vector parameter is larger than the IDQD, then the vector parameter space is not identifiable, and moreover, there are infinitely many nonidentifiable vector parameter points in the vector parameter space. Thus, the IDQD indicates a quantization induced fundamental limitation on the estimation capabilities of the quantized estimation system with regard to the identifiability condition. 

\item We next show that there is no general equivalence between the just described quantization induced FIM singularity and the quantization induced nonidentifiability of the vector parameter space. In particular, there exist some cases where the necessary conditions for the existence of the FIM  do not hold, but the quantization induced nonidentifiability of the vector parameter space can still be guaranteed. However, if the FIM exists, the condition that the dimension of the desired vector parameter is larger than the IDQD gives rise to both quantization induced issues.

\item Some further investigations into quantization induced nonidentifiability are carried out.  We show that in some cases where the dimension of the desired vector parameter is larger than the IDQD, every vector parameter point in the quantization induced nonidentifiable vector parameter space is nonidentifiable, while in some other cases, only some vector parameter points in the quantization induced nonidentifiable vector parameter space are identifiable. Thus the quantization induced FIM singularity does not necessarily determine the identifiability of the vector parameter point although it does determine the identifiability of the vector parameter space. Moreover, we show that the cardinality of a set of vector parameter points in the quantization induced nonidentifiable vector parameter space which are as likely as each other based on the observations can be as small as $1$ and can also be as large as uncountably infinite.

\item Finally, as opposed to our previous general results, we consider scenarios where some commonly assumed specific assumptions on the statistical models of the observations are made.  It is shown that under the assumptions, the fundamental limitation of the quantization system becomes more limiting. A smaller dimension of the vector parameter, called the refined IDQD (rIDQD), will guarantee the FIM singularity  and the nonidentifiability of the vector parameter space.  
\end{asparaenum}




\subsection{Related Work}

The information-regularity condition and the identifiability condition have been successfully applied in several engineering disciplines, including statistical inference, control theory, and array processing, see \cite{hero1996exploring, stoica1982non, hochwald1996identifiability, stoica2001parameter, paulino1994identifiability} for examples. 
Previous work has illuminated an intimate link between the nonsingularity of the FIM and the local identifiability of the desired parameter \cite{stoica1982non,  hochwald1996identifiability, stoica2001parameter, paulino1994identifiability,  hochwald1997identifiability, rothenberg1971identification,  bowden1973theory}. Local identifiability implies identifiability in an open neighborhood of the true value of the desired parameter, and is weaker than the identifiability discussed in this paper which is often called global identifiability. The author of \cite{rothenberg1971identification} shows that if the rank of the FIM is constant over some open neighborhood of the desired parameter, then the nonsingularity of the FIM is equivalent to the local identifiability of the desired parameter. For normal distributions, the work in \cite{hochwald1996identifiability} provides some other conditions which also guarantee the equivalence between the nonsingularity of the FIM and the local identifiability of the desired parameter. However,  it can be shown that quantization induced singularity of the FIM does not generally imply a lack of local identifiability.


The relationship between the identifiability and the dimension of the vector parameter to be estimated  has  been studied in the area of array processing for a particular class of multivariate Gaussian distributed signal models \cite{hochwald1996identifiability}. For the particular class of models considered in \cite{hochwald1996identifiability}, the task of examining the identifiability can be simplified to examining whether different values of the parameters give rise to different values of the covariance matrix. However, there are major differences between the work in \cite{hochwald1996identifiability} and that in this paper. First and foremost, the array processing models considered in \cite{hochwald1996identifiability} do not employ quantization which is the focus of our work.
Moreover, we don't make any assumption on the model of the received signals, and our results hold for arbitrary statistical models of the observations, arbitrary value of the desired vector parameter, and arbitrary quantizer designs.


\subsection{Notation and Organization}

In this paper, bold upper case letters and bold lower case letters are used to represent matrices and column vectors respectively.
The symbol ${\bf 1} $ stands for the all-one column vector, and ${\bf 0} $ for the all-zero column vector. 
For any set $\cal S$, $|{\cal S}|$ represents the number of elements in the set $\cal S$. For any given $L$, $\mathbbm{R}^L$ denotes the set of all $L$-tuples real numbers. The rank and expectation operators are denoted by ${{\rm{rank}(\cdot)}}$ and ${\mathbbm{E}}\left(  \cdot  \right)$ respectively.

The remainder of the paper is organized as follows. A general quantized  estimation system and its IDQD are introduced in \emph{Section \ref{Section_System_Model}}. In \emph{Section \ref{Section_Information_Regularity_Estimation_Capacity}}, the impact of  quantization on the information-regularity condition is investigated. \emph{Section \ref{Section_Nonidentifiability_Estimation_Capacity}} studies the impact of  quantization on the identifiability condition. The specialization of the results in \emph{Section \ref{Section_Information_Regularity_Estimation_Capacity}} and \emph{Section \ref{Section_Nonidentifiability_Estimation_Capacity}} to cases with some commonly assumed assumptions is considered in \emph{Section \ref{Section_Estimation_Capacity_Assumptions}}.
Finally, \emph{Section \ref{Section_Conclusion}} provides our conclusions.

\section{Quantized Parameter Estimation System Model and Inestimable Dimension for Quantized Data}
\label{Section_System_Model}

Consider an $N$-sensor system as shown in Fig. \ref{Fig_estimation_system} where the $j$-th sensor\footnote{It should be noted that if the data comes from something other than a sensor, the results still apply.} produces a $K$-dimensional vector ${\bf x}_{j}$. The statistical description of ${\bf x}_{j}$ depends on a $D_{{\boldsymbol{\theta }}}$-dimensional vector parameter ${\boldsymbol{\theta }} \in {\boldsymbol{\Theta}} \subset {\mathbbm{R}}^{D_{{\boldsymbol{\theta }}}}$ that we wish to estimate.   As a generalization\footnote{The case of performing $L_j$ scalar quantizations at each sensor, is common, as is pure vector quantization with $L_j=1$, for example.}  
to the standard quantized parameter estimation system, partition the observation vector into $L_j$ parts as 
\begin{equation}
{{\bf{x}}_j} = {\left[ {{\bf{x}}_{j1}^T,{\bf{x}}_{j2}^T,...,{\bf{x}}_{j{L_j}}^T} \right]^T}.
\end{equation}
Next the $l$-th observation subvector ${\bf x}_{jl}$ of ${\bf x}_{j}$ is converted to the quantized value $u_{jl}$ by employing the vector quantizer $\gamma_{jl}$ using 
\begin{equation} \label{u_jl}
{u_{jl}}  \buildrel \Delta \over = {\gamma _{jl}}\left( {{{\bf{x}}_{jl}}} \right) = \sum\limits_{r = 1}^{{R_{jl}}} {r \ {\mathbbm{1}}\left\{ {{{\bf{x}}_{jl}} \in I_{jl}^{(r)}} \right\}},
\end{equation}
where $ {\mathbbm{1}}\{ {{{\bf{x}}_{jl}} \in I_{jl}^{(r)}} \} = 1 $  if $ {{{\bf{x}}_{jl}} \in I_{jl}^{(r)}} $ and otherwise it is zero. 
Thus ${\gamma _{jl}}$ is an $R_{jl}$-level vector quantizer with given quantization regions $\{ {I_{jl}^{(r)}} \}_{r = 1}^{{R_{jl}}}$ which are disjoint and cover the domain of ${\gamma _{jl}}$. Next, we collect all the quantized data corresponding to ${\bf x}_j$, into ${\bf u}_j$ which we call the superquantized vector, such that  
\begin{align} \notag
{{\bf{u}}_j}  & = {\left[ {{u_{j1}},{u_{j2}},...,{u_{j{L_j}}}} \right]^T} \\ \notag
& \buildrel \Delta \over = {\Gamma _j}\left( {{{\bf{x}}_j}} \right) \\
& = {\left[ {{\gamma _{j1}}\left( {{{\bf{x}}_{j1}}} \right),{\gamma _{j2}}\left( {{{\bf{x}}_{j2}}} \right),...,{\gamma _{j{L_j}}}\left( {{{\bf{x}}_{j{L_j}}}} \right)} \right]^T}.
\end{align}
We assume for simplicity that the quantities $ {{\bf{u}}_j}, j=1,2,...,N $ are transmitted without error to the FC to be used for estimating ${\boldsymbol{\theta}} $. 


Without loss of generality, we can assume that the observation vectors $\{{\bf x}_{j}\}_{j=1}^{N}$ are independent, but the elements of ${\bf x}_{j}$ are not necessarily independent for each $j$. This is because the scenarios where $\{{\bf x}_{j}\}_{j=1}^{N}$ are not independent can be considered as a special case of the system  which only consists of $1$ sensor, since the observation vectors $\{{\bf x}_{j}\}_{j=1}^{N}$ at different sensors are just a partition of the overall observation vector ${[ {{\bf{x}}_1^{T},{\bf{x}}_2^T,...,{\bf{x}}_N^T} ]^T}$.  We assume that ${\bf x}_{j}$ follows a statistical model  $( {\mathscr{X}}_j, {{\mathscr F}_j}, {{{ \mathscr P}}}_j^{\boldsymbol{\theta }} )$ for each $j$, where ${\mathscr{X}}_j$ is some set endowed with a ${\sigma}$-algebra ${\mathscr F}_j$.  The probability measure ${{{ \mathscr P }}}_j^{\boldsymbol{\theta }}$ of ${\bf x}_{j}$ belongs to a family of probability measures $\{ {{\mathscr P}}_{j}^{\boldsymbol{\theta }}: {\boldsymbol{\theta }} \in {\boldsymbol{\Theta }}  \}$ on  $({\mathscr{X}}_j, {\mathscr F}_j)$ indexed by a $D_{{\boldsymbol{\theta }}}$-dimensional vector parameter ${\boldsymbol{\theta }}$ lying in a set ${\boldsymbol{\Theta}} \subset {\mathbbm{R}}^{D_{{\boldsymbol{\theta }}}}$. The superquantizer ${\Gamma _j}:({\mathscr{X}}_j, {\mathscr F}_j) \to ( {{{\mathbbm{R}}^{{L_j}}}, {\mathscr B}_j} )$ is a measurable function with respect to ${\mathscr F}_j$ and ${\mathscr B}_j$ for all $j$, where ${\mathscr B}_j$ is the Borel algebra on  ${\mathbbm{R}}^{{L_j}}$.

Before proceeding, we define a critical quantity, which is called the inestimable dimension for quantized data (IDQD) of the quantized estimation system.
\begin{definition}[Inestimable Dimension for Quantized Data] \label{Definition_estimation_capacity}
	The IDQD $\lambda \left( N,\left\{ {{R_{jl}}} \right\} \right)$ of the quantized estimation system described above is defined as
	\begin{equation} \label{Estimation_Capacity_Vector_Quantization}
	\lambda \left( {N,\left\{ {{R_{jl}}} \right\}} \right) \buildrel \Delta \over = \sum\limits_{j = 1}^N {\prod\limits_{l = 1}^{{L_j}} {{R_{jl}}} }  - N
	\end{equation}
\end{definition}

It is seen from (\ref{Estimation_Capacity_Vector_Quantization}) that the defined IDQD $\lambda \left( N,\left\{ {{R_{jl}}} \right\} \right)$ of the quantized estimation system does not depend on the value of the desired vector parameter $\boldsymbol{\theta}$, the statistical models $\{( {\mathscr{X}}_j, {{\mathscr F}_j}, {{{ \mathscr P}}}_j^{\boldsymbol{\theta }} )\}$ and the quantization regions $\{I_{jl}^{(r)}\}$, but is only determined by the number $N$ of sensors and the numbers $\{R_{jl}\}$ of the quantization levels employed at the sensors. In the following, we will show that there is a close link between the IDQD and the estimation capabilities of the quantized estimation system in terms of the information-regularity condition and the identifiability condition.


\section{IDQD and Information-Regularity Condition}
\label{Section_Information_Regularity_Estimation_Capacity}

In this section, we first formulate the FIM for estimating $\boldsymbol{\theta}$, and then based on the expression of the FIM, 
we show that the IDQD of the quantized estimation system describes a fundamental limitation of the quantized estimation system with respect to the information-regularity condition.

Let ${{\cal S}_j}$ denote the set of all possible outcomes of the $j$-th superquantizer ${\Gamma}_j$
\begin{equation} \label{S_p}
{{\cal S}_j} = \left\{ {{\bf{s}}_1^{(j)},{\bf{s}}_2^{(j)},...,{\bf{s}}_{|{{\cal S}_j}|}^{(j)}} \right\}.
\end{equation}
It is clear that the size of ${\cal S}_j$ can be written as 
\begin{equation} \label{size_Sp}
\left| {{{\cal S}_j}} \right| = \prod\limits_{l = 1}^{{L_j}} {{R_{jl}}}. 
\end{equation}
Let ${\bf{u}}$ denote a vector containing all the quantized data $\{{\bf u}_j\}$ received at the FC
\begin{equation} \label{u}
{\bf{u}} \buildrel \Delta \over = {\left[ {{\bf{u}}_1^T,{\bf{u}}_2^T,...,{\bf{u}}_N^T} \right]^T}.
\end{equation}  
For any given quantized data ${\bf{u}}$ received at the FC, the log-likelihood function $L\left( {\boldsymbol{\theta }} \right)$ can be expressed as\footnote{Note that if ${q}_j^{({\bf{s}})}\left( {\boldsymbol{\theta }} \right)=0$ for some $j$ and $\bf s$, then the corresponding summand in (\ref{LLR_vector}) should be eliminated in computing (\ref{LLR_vector}).} 
\begin{align} \notag
L\left( {\boldsymbol{\theta }} \right) & \buildrel \Delta \over = \ln \Pr \left( {{\bf{u}}\left| {\boldsymbol{\theta }} \right.} \right) \\ \notag
& = \ln \prod\limits_{j = 1}^N {\Pr \left( {{{\bf{u}}_j}\left| {\boldsymbol{\theta }} \right.} \right)} \\  \label{LLR_vector}
&  = \sum\limits_{j = 1}^N {\sum\limits_{{\bf{s}} \in {{\cal S}_j}} {{\mathbbm{1}}\left\{ {{{\bf{u}}_j} = {\bf{s}}} \right\}\ln q_j^{({\bf{s}})}\left( {\boldsymbol{\theta }} \right)} } 
\end{align}
where  $\forall j$, ${q}_j^{({\bf{s}})}\left( {\boldsymbol{\theta }} \right)$ is defined as
\begin{equation} \label{q_j_s}
q_j^{({\bf{s}})}\left( {\boldsymbol{\theta }} \right)  \buildrel \Delta \over = {\mathscr P}_j^{\boldsymbol{\theta }}\left( {{\Gamma _j}\left( {{{\bf{x}}_j}} \right) = {\bf{s}}} \right)
\end{equation}
for any given vector ${\bf{s}} \in {\cal S}_j$.

Define the following assumptions. 
\begin{assumption} \label{A1_Open_Set_Domain}
	For the set $\boldsymbol{\Theta}$ in ${\mathbbm{R}}^{D_{\boldsymbol{\theta}}}$, the interior of $\boldsymbol{\Theta}$ is not empty.
\end{assumption}
\begin{assumption} \label{A2_differentiable}
	For all $j$ and ${\bf s}$, ${q_{j}^{({\bf s})}\left( {\boldsymbol{\theta}} \right)}$ in (\ref{q_j_s}) is twice differentiable with respect to ${\boldsymbol{\theta}}$ for all ${\boldsymbol{\theta}} \in {\boldsymbol{\Theta}}$.
\end{assumption}
Note that \emph{Assumption \ref{A1_Open_Set_Domain}} and \emph{Assumption \ref{A2_differentiable}} are called regularity conditions and are commonly adopted in the signal processing literature \cite{poor1994introduction}. 

To gain insights into whether the information-regularity condition is satisfied, we first explore the FIM ${\bf{J}}(\boldsymbol{\theta})$ for estimating $\boldsymbol{\theta}$. Under \emph{Assumption \ref{A1_Open_Set_Domain}} and \emph{Assumption \ref{A2_differentiable}}, the $(l,m)$-th element of the FIM ${\bf{J}}(\boldsymbol{\theta})$ is defined as \cite{poor1994introduction}
\begin{equation}
{\left[ {{\bf{J}}\left( {\boldsymbol{\theta }} \right)} \right]_{l,m}} \buildrel \Delta \over =  - {\mathbbm{E}}\left\{ {\frac{{{\partial ^2}L\left( {\boldsymbol{\theta }} \right)}}{{\partial {\theta _l}\partial {\theta _m}}}} \right\}, 
\end{equation}
where $\theta_l$ and $\theta_m$ denote the $l$-th and $m$-th elements of $\boldsymbol{\theta}$ respectively. Hence, by employing (\ref{LLR_vector}), ${\bf{J}}(\boldsymbol{\theta})$ can be expressed as
\begin{equation} \label{FIM_vector_quantization}
{\bf{J}}\left( {\boldsymbol{\theta }} \right) = \sum\limits_{j = 1}^N {\sum\limits_{{\bf{s}} \in {{\cal S}_j}} {\frac{1}{{q_j^{({\bf{s}})}\left( {\boldsymbol{\theta }} \right)}}\frac{{\partial q_j^{({\bf{s}})}\left( {\boldsymbol{\theta }} \right)}}{{\partial {\boldsymbol{\theta }}}}{{\left[ {\frac{{\partial q_j^{({\bf{s}})}\left( {\boldsymbol{\theta }} \right)}}{{\partial {\boldsymbol{\theta }}}}} \right]}^T}} } .
\end{equation}

By employing (\ref{FIM_vector_quantization}), we can obtain the following theorem with regard to the singularity of the FIM.
\begin{theorem} \label{Theorem_Singularity_vector_quantization}
	Under \emph{Assumption \ref{A1_Open_Set_Domain}} and \emph{Assumption \ref{A2_differentiable}}, for any given $\boldsymbol{\theta}$,  any quantization regions $\{ {I_{jl}^{(r)}} \}$ and any statistical models $\{( {\mathscr{X}}_j, {{\mathscr F}_j}, {{{ \mathscr P}}}_j^{\boldsymbol{\theta }} )\}$, the FIM ${\bf{J}}\left(  {\boldsymbol{\theta}}   \right) $ described in (\ref{FIM_vector_quantization}) is singular, if the dimension $D_{\boldsymbol{\theta}}$ of the vector parameter $\boldsymbol{\theta}$ is greater than the IDQD, i.e., 
	\begin{equation} \label{Define_lambda_vector_quantization}
	D_{\boldsymbol{\theta}} > \lambda \left( {N,\left\{ {{R_{jl}}} \right\}} \right),
	\end{equation}
	where $\lambda ( {N,\{ {{R_{jl}}} \}} )$ is defined in (\ref{Estimation_Capacity_Vector_Quantization}).
\end{theorem}
\begin{IEEEproof}
	Refer to Appendix \ref{proof_Theorem_Singularity_vector_quantization}.
\end{IEEEproof}

\emph{Theorem \ref{Theorem_Singularity_vector_quantization}} reveals a fundamental limitation when utilizing quantized data for estimating a vector parameter, and sheds light on the preliminary design of a quantized estimation system. To be specific, the quantization and sensing approach employed should guarantee that the IDQD of the quantized estimation system, $\lambda \left( {N,\left\{ {{R_{jl}}} \right\}} \right)$,  is larger than or equal to the dimension of the vector parameter of interest. 
Otherwise, the FIM for estimating the vector parameter of interest is definitely singular for any $\boldsymbol{\theta}$, $\{I_{jl}^{(r)}\}$ and  $\{( {\mathscr{X}}_j, {{\mathscr F}_j}, {{{ \mathscr P}}}_j^{\boldsymbol{\theta }} )\}$. Since this fundamental limitation is quantization induced, we refer to this singularity of the FIM which is caused by the condition $D_{\boldsymbol{\theta}} > \lambda \left( {N,\left\{ {{R_{jl}}} \right\}} \right)$ as quantization induced singularity. 
In order to alleviate this undesirable outcome, it is seen from the definition of $\lambda \left( {N,\left\{ {{R_{jl}}} \right\}} \right)$ in (\ref{Estimation_Capacity_Vector_Quantization}) that one can employ finer quantizers $\{\gamma_{jl}\}$ with larger $\{R_{jl}\}$ or add more sensors in the system. However, enlarging $\{R_{jl}\}$ can significantly increase the data rate from each sensor to the FC, and adding more sensors increases the cost of quantized estimation system. 

If the dimension of the vector parameter of interest is smaller than the IDQD $\lambda \left( {N,\left\{ {{R_{jl}}} \right\}} \right)$, then it is possible that the FIM for estimating the vector parameter of interest is nonsingular for some $\boldsymbol{\theta}$, $\{ {I_{jl}^{(r)}} \}$ and  $\{( {\mathscr{X}}_j, {{\mathscr F}_j}, {{{ \mathscr P}}}_j^{\boldsymbol{\theta }} )\}$. However, in this case the singularity depends on $\boldsymbol{\theta}$, $\{ {I_{jl}^{(r)}} \}$ and  $\{( {\mathscr{X}}_j, {{\mathscr F}_j}, {{{ \mathscr P}}}_j^{\boldsymbol{\theta }} )\}$. Hence, the condition $D_{\boldsymbol{\theta}} > \lambda \left( {N,\left\{ {{R_{jl}}} \right\}} \right)$ is not generally necessary for guaranteeing the singularity of the FIM.

\section{IDQD and Identifiability Condition}
\label{Section_Nonidentifiability_Estimation_Capacity}

In this section, we study the relationship between the identifiability of the  vector parameter space $\boldsymbol{\Theta}$  and the IDQD of the quantized estimation system defined in (\ref{Estimation_Capacity_Vector_Quantization}). 
 The definitions of identifiability which are employed in this paper are first described. Then, we show that the IDQD of the quantized estimation system describes a fundamental limitation of the quantized estimation system with respect to the identifiability of the vector parameter space.




\subsection{IDQD and Identifiability of the Vector Parameter Space}
\label{Section_Nonidentifiable_Vector_Parameter_Space}

In order to characterize the impact of the quantization on the identifiability condition, we first formally give the following definitions with regard to the identifiability.

\begin{definition}[Observationally Equivalent \cite{rothenberg1971identification}] \label{Define_Observationally_Equivalent}
	Two distinct vector parameter points in ${\boldsymbol{\Theta}}$, ${\boldsymbol{\theta}}$ and  ${\boldsymbol{\theta}}'$, are said to be observationally equivalent if $\Pr \left( {{\bf{u}}\left| {\boldsymbol{\theta }} \right.} \right) = \Pr \left( {{\bf{u}}\left| {{{\boldsymbol{\theta }}}}' \right.} \right)$ for all possible $\bf{u}$ in (\ref{u}).
\end{definition}

\begin{definition}[Identifiable Vector Parameter Point \cite{rothenberg1971identification}] \label{Define_Identifiable_Parameter_Point}
	The vector parameter point ${\boldsymbol{\theta}} \in {\boldsymbol{\Theta}}$ is called identifiable,
	if there is no other ${\boldsymbol{\theta}}' \in {{\boldsymbol{\Theta }}\backslash \{{\boldsymbol{\theta }}\}}$ which is observationally equivalent to ${\boldsymbol{\theta}}$.
\end{definition}

\begin{definition}[Identifiable Vector Parameter Space] \label{Define_Identifiable_Parameter_Space}
	The vector parameter space $\boldsymbol{\Theta}$ is considered identifiable, 
	if every vector parameter point $\boldsymbol{\theta} \in {\boldsymbol{\Theta}}$ is identifiable.    
\end{definition}

It is worth pointing out that in some literature, if a parameter is said to be identifiable, it means that the parameter space ${\boldsymbol{\Theta}}$ is identifiable, for instance, see \cite{basu1983identifiability, lehmann1999elements}. In general, for a meaningful estimation problem, it is necessary that  the parameter space ${\boldsymbol{\Theta}}$ is identifiable.

Let ${\mathcal{A}} \buildrel \Delta \over = \left\{ {{{\bf{a}}_1},{{\bf{a}}_2},...,{{\bf{a}}_{{D_{\bf{u}}}}}} \right\}$ denote the set of all possible realizations of $\bf u$ in (\ref{u}), where the number of all possible realizations is
\begin{equation} \label{D_u}
{D_{\bf{u}}} = \prod\limits_{j = 1}^N {\prod\limits_{l = 1}^{{L_j}} {{R_{jl}}} }. 
\end{equation} 

According to \emph{Definition \ref{Define_Observationally_Equivalent}} and  \emph{Definition \ref{Define_Identifiable_Parameter_Space}}, the vector parameter space $\boldsymbol{\Theta}$ is identifiable if and only if there are no distinct vector parameter points ${\boldsymbol{\theta}}_1$ and ${\boldsymbol{\theta}}_2$ in ${\boldsymbol{\Theta}}$ such that $\Pr \left( {{\bf{u}}\left| {\boldsymbol{\theta }}_1 \right.} \right) = \Pr \left( {{\bf{u}}\left| {{{\boldsymbol{\theta }}_2}} \right.} \right)$ for all  $\bf{u} \in {\mathcal{A}}$. In other words,  the vector parameter space $\boldsymbol{\Theta}$ is identifiable if and only if the mapping
\begin{equation} \label{phi_u}
\begin{aligned} 
{\varphi _{\bf{u}}}:  {\boldsymbol{\Theta }} & \longrightarrow  {{\mathbbm{R}}^{{D_{\bf{u}}}}}\\ 
{\boldsymbol{\theta }}  & \longmapsto  {\left[ {\Pr \left( {{{\bf{a}}_1}\left| {\boldsymbol{\theta }} \right.} \right),\Pr \left( {{{\bf{a}}_2}\left| {\boldsymbol{\theta }} \right.} \right),...,\Pr \left( {{{\bf{a}}_{{D_{\bf{u}}}}}\left| {\boldsymbol{\theta }} \right.} \right)} \right]^T}
\end{aligned} 
\end{equation}
is injective. Thus, we can examine the injectivity of the mapping ${\varphi _{\bf{u}}}$ in (\ref{phi_u}) to investigate the identifiability condition. On the other hand, the following lemmas simplify this investigation.

Let us define a $(\sum\nolimits_{j = 1}^N {\prod\nolimits_{l = 1}^{{L_j}} {{R_{jl}}} }  - N)$-dimensional vector ${\boldsymbol{\Psi }}( {\boldsymbol{\theta }} )$ 
\begin{equation} \label{Psi}
{\boldsymbol{\Psi }}\left( {\boldsymbol{\theta }} \right) \buildrel \Delta \over = {\left[ {{{\boldsymbol{\psi }}_1}{{\left( {\boldsymbol{\theta }} \right)}^T},{{\boldsymbol{\psi }}_2}{{\left( {\boldsymbol{\theta }} \right)}^T},...,{{\boldsymbol{\psi }}_N}{{\left( {\boldsymbol{\theta }} \right)}^T}} \right]^T},
\end{equation}
where for each $j$, ${{\boldsymbol{\psi }}_j}\left( {\boldsymbol{\theta }} \right)$ is defined as
\begin{equation} \label{psi}
{{\boldsymbol{\psi }}_j}\left( {\boldsymbol{\theta }} \right) \buildrel \Delta \over = {\left[ {q_j^{({\bf{s}}_1^{(j)})}\left( {\boldsymbol{\theta }} \right),q_j^{({\bf{s}}_2^{(j)})}\left( {\boldsymbol{\theta }} \right),...,q_j^{({\bf{s}}_{|{{\cal S}_j}| - 1}^{(j)})}\left( {\boldsymbol{\theta }} \right)} \right]^T},
\end{equation}
$q_{j}^{({\bf s})}\left( {\boldsymbol{\theta }} \right)$ is defined in (\ref{q_j_s}), and ${\bf s}_i^{(j)}$ is defined in (\ref{S_p}) for all $i=1,2,...,|{{\cal S}_j}| - 1$.

\begin{lemma} \label{Lemma_iff_identifiability_condition}
	The mapping ${{{\varphi}_{\bf u}} }$ in (\ref{phi_u}) is injective if and only if the mapping
	\begin{equation} \label{map_reduce_dimension_Psi}
	\begin{aligned} 
	{\boldsymbol{\Psi }} :  {\boldsymbol{\Theta }} & \longrightarrow  {{\mathbbm{R}}^{\sum\limits_{j = 1}^N {\prod\limits_{l = 1}^{{L_j}} {{R_{jl}}} }  - N}}\\ 
	{\boldsymbol{\theta}} &  \longmapsto  {\boldsymbol{\Psi }}\left( {\boldsymbol{\theta }} \right)
	\end{aligned}
	\end{equation}
	is injective. Therefore, a necessary and sufficient condition under which the vector parameter space ${\boldsymbol{\Theta}}$ is identifiable is that the mapping ${\boldsymbol{\Psi }}$ in (\ref{map_reduce_dimension_Psi}) is injective. Moreover, the dimension of the vector ${\boldsymbol{\Psi }}\left( {\boldsymbol{\theta }} \right)$ in (\ref{Psi}) is strictly smaller than that of ${{{\varphi}_{\bf u}} }\left( {\boldsymbol{\theta }} \right)$ in (\ref{phi_u}) for any given $N$ and $\{R_{jl}\}$.
\end{lemma}
\begin{IEEEproof}
	Refer to Appendix \ref{proof_Lemma_iff_identifiability_condition}.
\end{IEEEproof}

As \emph{Lemma \ref{Lemma_iff_identifiability_condition}} demonstrates, the identifiability of the vector parameter space $\boldsymbol{\Theta}$ can also be determined by the injectivity of the mapping ${\boldsymbol{\Psi }}$ in (\ref{map_reduce_dimension_Psi}). To this end, we only need to inspect the injectivity of the mapping ${\boldsymbol{\Psi }}$ in (\ref{map_reduce_dimension_Psi}) to investigate the identifiability of the vector parameter space. What's more, it is seen that the dimension of the vector ${\boldsymbol{\Psi }}\left( {\boldsymbol{\theta }} \right)$ in (\ref{Psi}) is precisely the IDQD of the quantized estimation system $\lambda ( {N,\{ {{R_{jl}}} \}} )$, which is shown to be strictly smaller than the dimension of ${{{\varphi}_{\bf u}} }\left( {\boldsymbol{\theta }} \right)$ in (\ref{phi_u}) for any given $N$ and $\{R_{jl}\}$.
In the following, we will show that because of the smaller dimension of ${\boldsymbol{\Psi }}\left( {\boldsymbol{\theta }} \right)$, inspecting the injectivity of the mapping ${\boldsymbol{\Psi }}$ in (\ref{map_reduce_dimension_Psi}) is easier than inspecting the injectivity of the mapping ${{{\varphi}_{\bf u}} }$ in (\ref{phi_u})
under the condition that $D_{\boldsymbol{\theta}} > \lambda \left( {N,\left\{ {{R_{jl}}} \right\}} \right)$.

Before proceeding, let us first introduce a helpful result in algebraic topology.
\begin{lemma}[Invariance of Domain Theorem \cite{tom2008algebraic}] \label{Lemma_domain_invariance}
	If $\cal U$ is an open subset of ${\mathbbm{R}}^n$ and $f: {\cal U} \to {\mathbbm{R}}^n$ is an injective continuous mapping, then ${\cal V}  \buildrel \Delta \over =  f({\cal U})$ is open in ${\mathbbm{R}}^n$, and $f$ is a homeomorphism between $\cal U$ and $\cal V$.
\end{lemma}

The proof of \emph{Lemma \ref{Lemma_domain_invariance}} can be found in \cite{tom2008algebraic}.
Next, we make the following assumption throughout this section, which is weaker than \emph{Assumption \ref{A2_differentiable}}.
\begin{assumption} \label{A3_continuity}
	For all $j$ and all ${\bf s}$, ${q_{j}^{({\bf s})}\left( {\boldsymbol{\theta}} \right)}$ in (\ref{q_j_s}) is continuous with respect to ${\boldsymbol{\theta}}$.
\end{assumption}

By employing \emph{Lemma \ref{Lemma_iff_identifiability_condition}} and \emph{Lemma \ref{Lemma_domain_invariance}}, we provide the following theorem with regard to the identifiability of the vector parameter space.

\begin{theorem} \label{Theorem_nonidentifiable_vector_parameter_space}
	Under \emph{Assumption \ref{A1_Open_Set_Domain}} and \emph{Assumption \ref{A3_continuity}}, for any given quantization regions $\{ {I_{jl}^{(r)}} \}$ and statistical models $\{( {\mathscr{X}}_j, {{\mathscr F}_j}, {{{ \mathscr P}}}_j^{\boldsymbol{\theta }} )\}$, if the dimension of the desired vector parameter $\boldsymbol{\theta}$ is larger than the IDQD of the quantized estimation system, i.e.,
	\begin{equation} \label{D_larger_than_estimation_capacity}
	D_{\boldsymbol{\theta}} > \lambda \left( {N,\left\{ {{R_{jl}}} \right\}} \right) 
	\end{equation} 
	then the vector parameter space $\boldsymbol{\Theta}$ is not identifiable. Moreover, for any open subset ${\cal U} \subset {\boldsymbol{\Theta}}$ in ${\mathbbm{R}}^{D_{\boldsymbol{\theta}}}$, there are infinitely many vector parameter points in ${\cal U}$ which are not identifiable.
\end{theorem}

\begin{IEEEproof}
	Refer to Appendix \ref{proof_Theorem_nonidentifiable_vector_parameter_space}
\end{IEEEproof}

\emph{Theorem \ref{Theorem_nonidentifiable_vector_parameter_space}} demonstrates that under \emph{Assumption \ref{A1_Open_Set_Domain}} and \emph{Assumption \ref{A3_continuity}}, for any given $\{ {I_{jl}^{(r)}} \}$ and $\{( {\mathscr{X}}_j, {{\mathscr F}_j}, {{{ \mathscr P}}}_j^{\boldsymbol{\theta }} )\}$, the condition $D_{\boldsymbol{\theta}} > \lambda \left( {N,\left\{ {{R_{jl}}} \right\}} \right)$ is sufficient to guarantee that there exist
infinitely many nonidentifiable vector parameter points in the vector parameter space ${\boldsymbol{\Theta}}$. Hence, the vector parameter space $\boldsymbol{\Theta}$ is not identifiable. This nonidentifiability of the vector parameter space is also quantization induced,
and doesn't depend on the statistical models of the observations and the design of the quantization regions. We refer to the nonidentifiability of the vector parameter space which is caused by the condition $D_{\boldsymbol{\theta}} > \lambda \left( {N,\left\{ {{R_{jl}}} \right\}} \right)$ as quantization induced nonidentifiability. However, it is worth mentioning that the condition $D_{\boldsymbol{\theta}} \le \lambda \left( {N,\left\{ {{R_{jl}}} \right\}} \right)$ cannot guarantee the identifiability of the vector parameter space, which is determined by the vector parameter space $\boldsymbol{\Theta}$, the quantizer designs $\{I_{jl}^{(r)}\}$, and the statistical models $\{( {\mathscr{X}}_j, {{\mathscr F}_j}, {{{ \mathscr P}}}_j^{\boldsymbol{\theta }} )\}$.

\subsection{Remarks on Quantization Induced Nonidentifiable Vector Parameter Space}

A particular note of interest is that \emph{Assumption \ref{A3_continuity}} employed in \emph{Theorem \ref{Theorem_nonidentifiable_vector_parameter_space}} is much weaker than \emph{Assumption \ref{A2_differentiable}} employed in \emph{Theorem \ref{Theorem_Singularity_vector_quantization}}.
The continuity of ${q_{j}^{({\bf s})}\left( {\boldsymbol{\theta}} \right)}$ assumed in \emph{Assumption \ref{A3_continuity}} is not enough to guarantee the existence of the FIM. Thus, in some cases where
the FIM for estimating the desired vector parameter doesn't exist, the quantization induced nonidentifiability of the vector parameter space can still be guaranteed by \emph{Theorem \ref{Theorem_nonidentifiable_vector_parameter_space}} under the condition that $D_{\boldsymbol{\theta}} > \lambda \left( {N,\left\{ {{R_{jl}}} \right\}} \right)$. Hence, in general, there is no equivalence between the quantization induced singularity of the FIM and the quantization induced nonidentifiability of the vector parameter space. However, if both \emph{Assumption \ref{A1_Open_Set_Domain}} and \emph{Assumption \ref{A2_differentiable}} hold, the condition that $D_{\boldsymbol{\theta}} > \lambda(N,\{R_{jl}\})$ ensures both singularity of the FIM and nonidentifiability of the vector parameter space.

According to \emph{Theorem \ref{Theorem_nonidentifiable_vector_parameter_space}}, we know that under \emph{Assumption \ref{A1_Open_Set_Domain}} and \emph{Assumption \ref{A3_continuity}}, for any given $\{ {I_{jl}^{(r)}} \}$ and $\{( {\mathscr{X}}_j, {{\mathscr F}_j}, {{{ \mathscr P}}}_j^{\boldsymbol{\theta }} )\}$, if the dimension of the desired vector parameter $\boldsymbol{\theta}$ is larger than the IDQD of the quantized estimation system,  
then there exist infinitely many nonidentifiable vector parameter points in ${\boldsymbol{\Theta}}$. 
However, there still remain two interesting questions which are not answered by \emph{Theorem \ref{Theorem_nonidentifiable_vector_parameter_space}}. First, although \emph{Theorem \ref{Theorem_nonidentifiable_vector_parameter_space}} shows that under the condition $D_{\boldsymbol{\theta}} > \lambda(N,\{R_{jl}\})$, there are infinitely many nonidentifiable vector parameter points in $\boldsymbol{\Theta}$, whether every vector parameter point in $\boldsymbol{\Theta}$ is nonidentifiable or not is still unknown. Second, under the condition $D_{\boldsymbol{\theta}} > \lambda(N,\{R_{jl}\})$, what is the cardinality of a given set of observationally equivalent points in the nonidentifiable vector parameter space?

In this subsection, we employ examples to show that in some cases, every vector parameter point in the vector parameter space is not identifiable, while in other cases, there exist some vector parameter points which are identifiable.  Moreover, the examples illustrate that under the condition $D_{\boldsymbol{\theta}} > \lambda(N,\{R_{jl}\})$, the cardinality of a set of observationally equivalent points can be very different for various cases.  It can be as small as $1$ and can also be as large as uncountably infinite.

\subsubsection{Every Vector Parameter Point in the Nonidentifiable Vector Parameter Space is Nonidentifiable and Every Set of Observationally Equivalent Points is Uncountable}
\ 

{\it Example 1}: Consider a quantized estimation system with $N=1$, $L_j=1$, and the dimension of ${\bf x}_{jl} $ is $1$ for all $l$.  In this case the single sensor makes a scalar observation which we denote as $x$ for simplicity.  The Gaussian assumed probability density function (pdf) of $x$ is
\begin{equation} \label{example_pdf}
f\left( {x\left| {\boldsymbol{\theta}}   \right.} \right) = \frac{1}{{\sqrt {2\pi {\beta }} }}{e^{ - \frac{{{{(x - \alpha )}^2}}}{{2{\beta }}}}},
\end{equation}
where the unknown vector parameter is ${\boldsymbol{\theta}} \buildrel \Delta \over =  [ \alpha, \beta]^T$. The vector parameter space is
\begin{equation} \label{Theta_domain_example}
{\boldsymbol{\Theta}}  = \left\{ {{{\left[ {\alpha ,{\beta }} \right]^T}}: \; {\alpha  \in {\mathbbm{R}}, \; {\beta \ge 0} }} \right\}.
\end{equation}
It is clear that the interior of ${\boldsymbol{\Theta}} $ is not empty.
We assume that the sensor employs a binary quantizer to convert $x$ to $u \in \{1,2\}$ by using
the nonempty quantization regions
\begin{equation} \label{everypoint_nonidentifiable_quantizer}
I^{(1)} = (a , b)  \text{ and } I^{(2)} = {\mathbbm{R}} \backslash I^{(1)},
\end{equation} 
for some $a$ and $b$ with $-\infty \le a < b \le \infty$.  This quantizer model is common and widely considered in recent literature, for instance, see  \cite{xiao2006distributed2, venkitasubramaniam2007quantization}.  
By the definition of the IDQD in (\ref{Estimation_Capacity_Vector_Quantization}), we can obtain 
\begin{equation}
\lambda \left( {N,\left\{ {{R_{jl}}} \right\}} \right) = \lambda \left( { 1 ,\left\{ 2 \right\}} \right) = 1 < 2 = {D_{\boldsymbol{\theta }}},
\end{equation}
and hence, by \emph{Theorem \ref{Theorem_nonidentifiable_vector_parameter_space}}, the vector parameter space ${\boldsymbol{\Theta}}$ is not identifiable.


\begin{proposition}  \label{Proposition_everypoint_nonidentifiable}
	For the quantized estimation system described in {\it Example 1}, if the sensor doesn't employ quantization, then every vector parameter point ${\boldsymbol{\theta}} \in {\boldsymbol{\Theta}}$ is identifiable. Since the vector parameter space $\boldsymbol{\Theta}$ is not identifiable when the quantizer is employed, the nonidentifiability of the vector parameter space is indeed quantization induced. Furthermore, if the sensor employs the quantizer with any given $I^{(1)}$ and $I^{(2)}$ with the forms in (\ref{everypoint_nonidentifiable_quantizer}), every vector parameter point in ${\boldsymbol{\Theta}}$ is not identifiable, and moreover, for any vector parameter point ${\boldsymbol{\theta}} \in {\boldsymbol{\Theta}}$, the set of vector parameter points which are observationally equivalent to ${\boldsymbol{\theta}}$ is uncountable.
\end{proposition}
\begin{IEEEproof}
	Refer to Appendix \ref{proof_Proposition_everypoint_nonidentifiable}.
\end{IEEEproof}

As \emph{Proposition \ref{Proposition_everypoint_nonidentifiable}} demonstrates, under the condition $D_{\boldsymbol{\theta}} > \lambda(N,\{R_{jl}\})$, there exist some cases where for any  vector parameter point ${\boldsymbol{\theta}} \in {\boldsymbol{\Theta}}$, the set of vector parameter points which are observationally equivalent to ${\boldsymbol{\theta}}$ is uncountable, and hence, every vector parameter point in ${\boldsymbol{\Theta}}$ is not identifiable.

\begin{figure}[htb]
	\vspace{0.028in}
	\centerline{
		\includegraphics[width=0.46\textwidth]{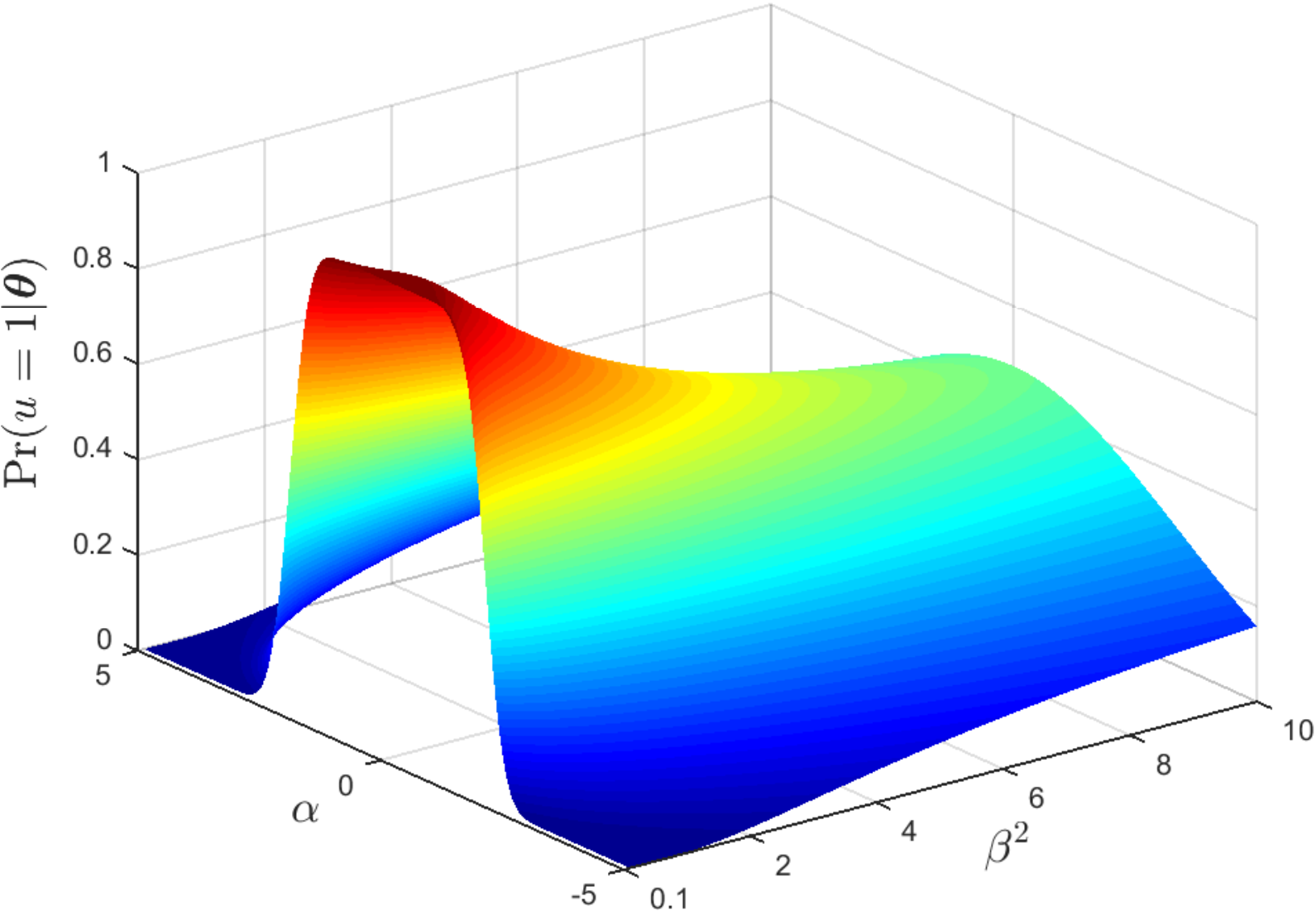}
	}
	\vspace{0.15in}
	\caption{$\Pr \left( {u = 1\left| {\boldsymbol{\theta }} \right.} \right)$ versus $\boldsymbol{\theta}$ for the case where $I^{(1)}=(-2,2)$ and $I^{(2)} = {\mathbbm{R}} \backslash I^{(1)}$.}
	\label{Fig_everypoint_nonidentifiability_surf}
\end{figure}

\begin{figure}[htb]
	\vspace{0.028in}
	\centerline{
		\includegraphics[width=0.46\textwidth]{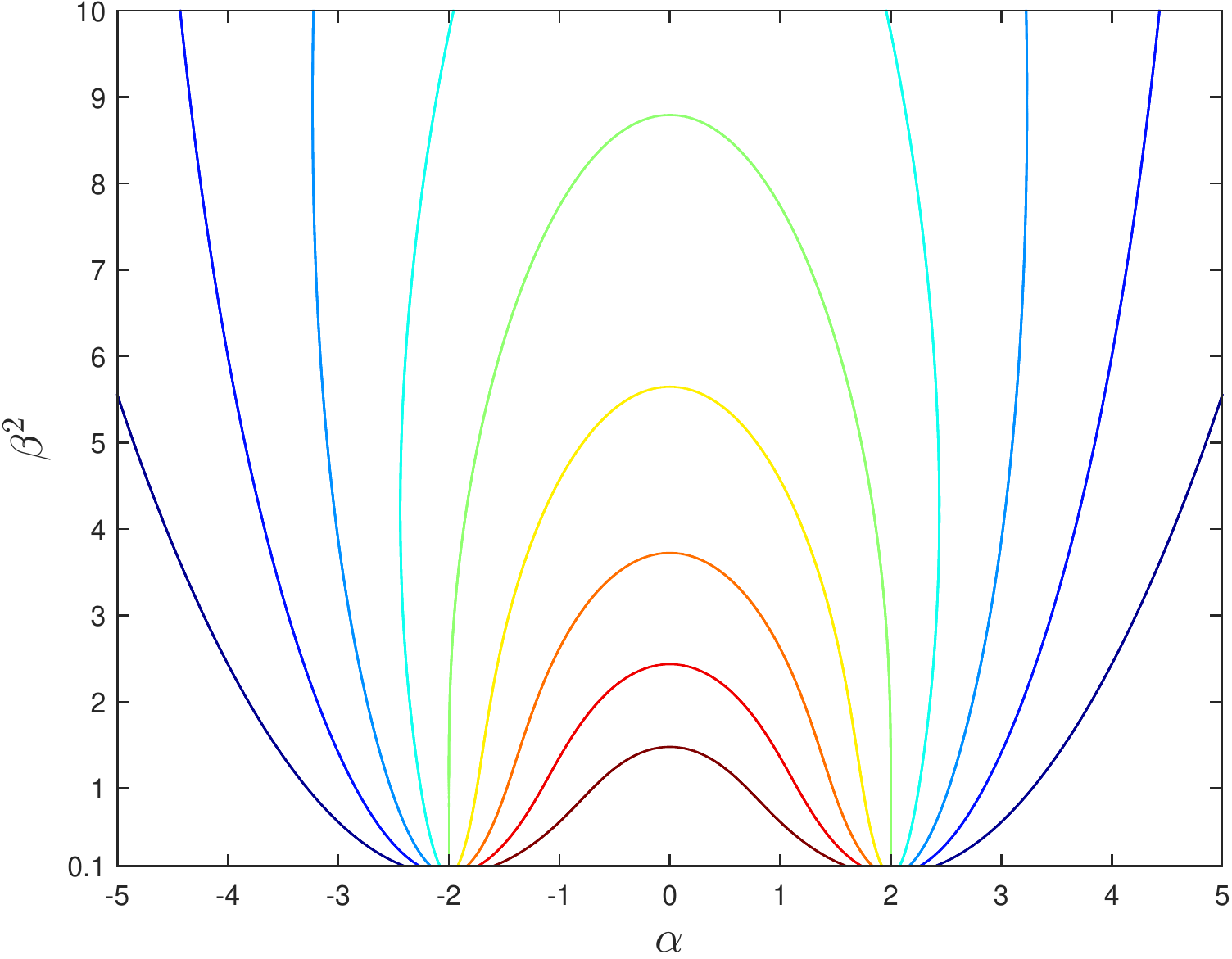}
	}
	\caption{Contours of $\Pr \left( {u = 1\left| {\boldsymbol{\theta }} \right.} \right)$ for the case where $I^{(1)}=(-2,2)$ and $I^{(2)} = {\mathbbm{R}} \backslash I^{(1)}$.}
	\label{Fig_everypoint_nonidentifiability_contour}
\end{figure}

To corroborate the theoretic analysis, we present some numerical results which illustrate the identifiability of the vector parameter points in ${\boldsymbol{\Theta}}$.
Fig. \ref{Fig_everypoint_nonidentifiability_surf} depicts the relationship between $\Pr \left( {u = 1\left| {\boldsymbol{\theta }} \right.} \right)$ and $\boldsymbol{\theta}$ for a particular case where $I^{(1)}=(-2,2)$ and $I^{(2)} = {\mathbbm{R}} \backslash I^{(1)}$, and Fig. \ref{Fig_everypoint_nonidentifiability_contour} shows the contour of  $\Pr \left( {u = 1\left| {\boldsymbol{\theta }} \right.} \right)$ for the same case. 
Since $\Pr \left( {u = 2\left| {\boldsymbol{\theta }} \right.} \right) = 1- \Pr \left( {u = 1\left| {\boldsymbol{\theta }} \right.} \right) $, we know that for a given ${\boldsymbol{\theta}}$ in $\boldsymbol{\Theta}$, if $\Pr \left( {u = 1\left| {\boldsymbol{\theta }}' \right.} \right) = \Pr \left( {u = 1\left| {\boldsymbol{\theta }} \right.} \right)$ for some other ${{\boldsymbol{\theta }}}'$ in $\boldsymbol{\Theta} \backslash \{\boldsymbol{\theta}\}$, then by \emph{Definition \ref{Define_Observationally_Equivalent}}, ${\boldsymbol{\theta}}'$ is observationally equivalent to ${{\boldsymbol{\theta }}}$, and hence, ${{\boldsymbol{\theta }}}$ is not identifiable. Therefore, every contour curve in Fig. \ref{Fig_everypoint_nonidentifiability_contour} illustrates a set of observationally equivalent vector parameter points. Moreover, it is easy to see from Fig. \ref{Fig_everypoint_nonidentifiability_surf} that every vector parameter point is not identifiable. 

\subsubsection{Existence of Identifiable Vector Parameter Point in the Nonidentifiable Vector Parameter Space}
\

{\it Example 2:} Consider a quantized estimation system with $N=1$ and $K=2$, where the observation ${\bf x} = [x_1, x_2]^T$ follows the distribution\footnote{ ${\mathcal{N}}({\boldsymbol{\theta}}, {\bf I})$ denotes a multivariate Gaussian distribution
	with mean vector ${\boldsymbol{\theta}}$ and covariance matrix ${\bf I}$, where ${\bf I}$ is the $2$-by-$2$ identity matrix.} ${\mathcal{N}}({\boldsymbol{\theta}}, {\bf I})$ with unknown vector parameter ${\boldsymbol{\theta}} \buildrel \Delta \over =  [ \theta_1, \theta_2]^T$. The vector parameter space ${\boldsymbol{\Theta}}$ is ${\mathbbm{R}}^2$ which is open. We assume that the sensor employs a binary vector quantizer to convert ${\bf x}$ to $u \in \{1,2\}$ by using the nonempty quantization regions
\begin{equation} \label{some_identifiable_point_quantizer}
I^{(1)} = (a_1 , b_1) \times (a_2, b_2)  \text{ and } I^{(2)} = {\mathbbm{R}}^2 \backslash I^{(1)},
\end{equation} 
for some $a_1$, $a_2$, $b_1$ and $b_2$, where $- \infty \le a_1 < b_1 \le \infty$ and $- \infty \le a_2 < b_2 \le \infty$. 

From the definition of the IDQD in (\ref{Estimation_Capacity_Vector_Quantization}), we can obtain 
\begin{equation}
\lambda \left( {N,\left\{ {{R_{jl}}} \right\}} \right) = \lambda \left( { 1 ,\left\{ 2 \right\}} \right) = 1 < 2 = {D_{\boldsymbol{\theta }}},
\end{equation}
and hence, by \emph{Theorem \ref{Theorem_nonidentifiable_vector_parameter_space}}, the vector parameter space ${\boldsymbol{\Theta}}$ is not identifiable.

\begin{proposition}  \label{Proposition_somepoint_identifiable}
	For the quantized estimation system described in {\it Example 2} and for any given $I^{(1)}$ and $I^{(2)}$ with the forms in (\ref{some_identifiable_point_quantizer}),
	 there exists an identifiable vector parameter point in ${\boldsymbol{\Theta}}$.
\end{proposition}
\begin{IEEEproof}
	Refer to Appendix \ref{proof_Proposition_somepoint_identifiable}.
\end{IEEEproof}

As illustrated by \emph{Proposition \ref{Proposition_somepoint_identifiable}}, under the condition $D_{\boldsymbol{\theta}} > \lambda(N,\{R_{jl}\})$, there exist some cases where the quantization induced nonidentifiable vector parameter space ${\boldsymbol{\Theta}}$ contains some identifiable vector parameter points. According to \emph{Definition \ref{Define_Identifiable_Parameter_Point}}, for any identifiable vector parameter point $\boldsymbol{\theta}$,  the set of points which are observationally equivalent to $\boldsymbol{\theta}$ consists of only one point, that is, $\boldsymbol{\theta}$.
Thus, under the condition $D_{\boldsymbol{\theta}} > \lambda(N,\{R_{jl}\})$, the cardinality of some observationally equivalent set in nonidentifiable vector parameter space ${\boldsymbol{\Theta}}$ can be as small as $1$ in some cases, since it is possible that some  vector parameter point in ${\boldsymbol{\Theta}}$ is identifiable. 
Furthermore, it is worth mentioning that under the condition $D_{\boldsymbol{\theta}} > \lambda \left( {N,\left\{ {{R_{jl}}} \right\}} \right) $, the FIM evaluated at any vector parameter point in ${\boldsymbol{\Theta}}$ is singular for any case. Hence, as \emph{Proposition \ref{Proposition_somepoint_identifiable}} demonstrates, the singularity of the FIM does not necessarily contradict the identifiability of the vector parameter point.

\begin{figure}[htb]
	\centerline{
		\includegraphics[width=0.46\textwidth]{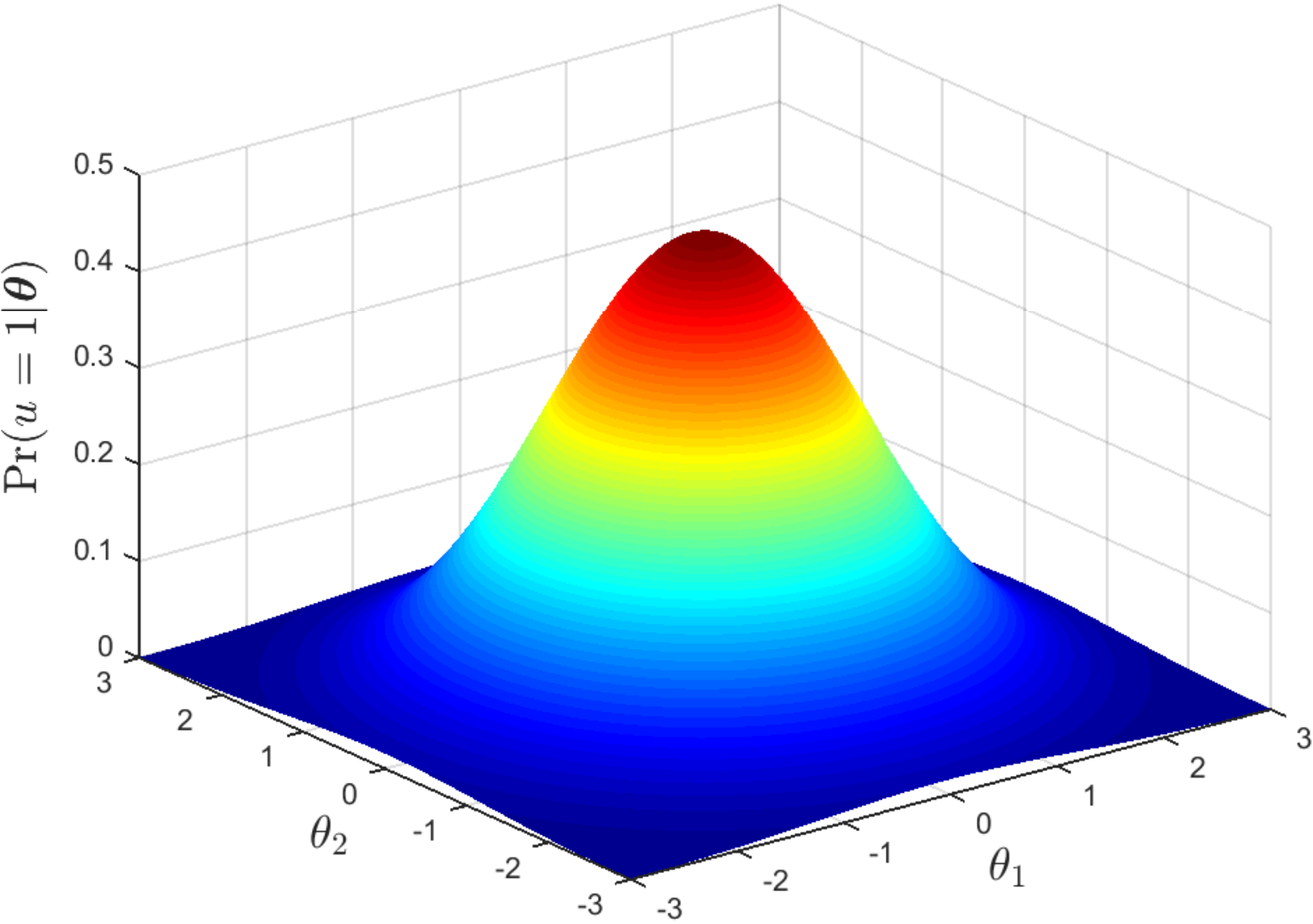}
	}
	\vspace{0.15in}
	\caption{$\Pr \left( {u = 1\left| {\boldsymbol{\theta }} \right.} \right)$ versus $\boldsymbol{\theta}$ for the case where  $a_1=a_2=-1$ and $b_1=b_2=1$.}
	\label{Fig_someypoint_identifiability}
\end{figure}

\begin{figure}[htb]
	\centerline{
		\includegraphics[width=0.46\textwidth]{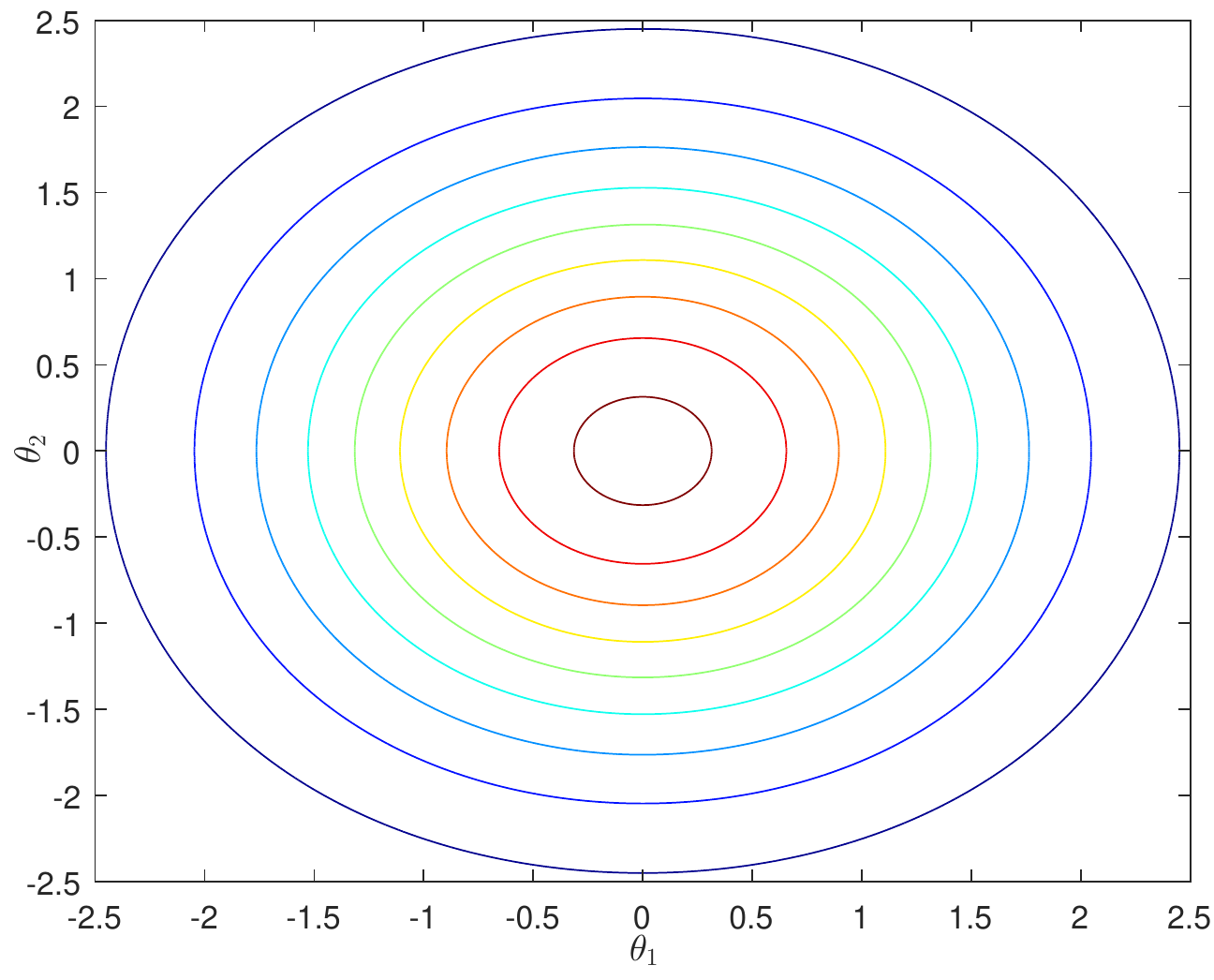}
	}
	\caption{Contours of $\Pr \left( {u = 1\left| {\boldsymbol{\theta }} \right.} \right)$ for the case where $a_1=a_2=-1$ and $b_1=b_2=1$.}
	\label{Fig_someypoint_identifiability_contour}
\end{figure}

Some numerical results for a particular case where $a_1=a_2=-1$ and $b_1=b_2=1$ are provided in support of the theoretical analysis.
Fig. \ref{Fig_someypoint_identifiability} depicts the value of  $\Pr \left( {u = 1\left| {\boldsymbol{\theta }} \right.} \right)$ for each vector parameter point in $\boldsymbol{\Theta}$, and Fig. \ref{Fig_someypoint_identifiability_contour} illustrates the contour of $\Pr \left( {u = 1\left| {\boldsymbol{\theta }} \right.} \right)$. Since $\Pr \left( {u = 2\left| {\boldsymbol{\theta }} \right.} \right) = 1- \Pr \left( {u = 1\left| {\boldsymbol{\theta }} \right.} \right) $, according to \emph{Definition \ref{Define_Identifiable_Parameter_Point}}, it is clear that if ${\boldsymbol{\theta}}$ is not identifiable, then there exists some other ${{\boldsymbol{\theta }}}'$ in $\boldsymbol{\Theta} \backslash \{\boldsymbol{\theta}\}$ such that $\Pr \left( {u = 1\left| {\boldsymbol{\theta }} \right.} \right) = \Pr \left( {u = 1\left| {\boldsymbol{\theta }}' \right.} \right)$. Fig. \ref{Fig_someypoint_identifiability} shows that $\Pr \left( {u = 1\left| {\boldsymbol{\theta }} \right.} \right)$ achieves its unique global maximum at ${\boldsymbol{\theta}} = {\bf 0}$, which demonstrates that ${\boldsymbol{\theta}} = {\bf 0}$ is an identifiable vector parameter point in this particular case.
It is seen from Fig. \ref{Fig_someypoint_identifiability_contour} that every set of observationally equivalent points forms a circle with the center at ${\boldsymbol{\theta}} = {\bf 0}$. Hence, except the set of points which are observationally equivalent to ${\boldsymbol{\theta}} = {\bf 0}$,  every other set of observationally equivalent points is uncountable.  

\section{IDQD with Additional Assumptions}
\label{Section_Estimation_Capacity_Assumptions}

In \emph{Section \ref{Section_Information_Regularity_Estimation_Capacity}} and \emph{Section \ref{Section_Nonidentifiability_Estimation_Capacity}}, we make no assumptions about the quantizers and the statistical model of the observations at each sensor. 
Hence, \emph{Theorem \ref{Theorem_Singularity_vector_quantization}} and \emph{Theorem \ref{Theorem_nonidentifiable_vector_parameter_space}} apply to any case with any $\{{\gamma}_{jl}\}$ and  $\{( {\mathscr{X}}_j, {{\mathscr F}_j}, {{{ \mathscr P}}}_j^{\boldsymbol{\theta }} )\}$. However, in general, the sufficient condition $D_{\boldsymbol{\theta}} > \lambda ( {N,\{ {{R_{jl}}} \}} )$ employed in \emph{Theorem \ref{Theorem_Singularity_vector_quantization}} and \emph{Theorem \ref{Theorem_nonidentifiable_vector_parameter_space}} for guaranteeing the FIM singularity  and the nonidentifiability of the vector parameter space is not strictly necessary. Hence, weaker conditions which still imply singularity and nonidentifiability when some additional assumptions are valid are of interest. In fact, under some common assumptions, we will show we can obtain a smaller IDQD compared to the results given in \emph{Theorem \ref{Theorem_Singularity_vector_quantization}} and \emph{Theorem \ref{Theorem_nonidentifiable_vector_parameter_space}} which better describes the limitations imposed by quantization under these assumptions. 

\subsection{Existence of Identical Sensor Observation Statistical Models and Identical Superquantizers}
\label{Section_Assumption_Identical_Model}

In this subsection, we consider the scenarios where the following assumption is valid.
\begin{assumption} \label{Assumption_identical_model}
	The statistical models of the observation vectors at some different sensors are known to be the same for all ${\boldsymbol{\theta }}$ such that 
the number of different statistical models at all the $N$ sensors is $P < N$. 
\end{assumption} 

Collect all the sensors indices that employ the $p$-th statistical model in the group ${\mathcal G}_p$ such that 
\begin{equation}
\{1,2,...,N\} = \mathop  \cup \limits_{p = 1}^P {{\cal G}_p}, \text{ and } {{\cal G}_p} \cap {{\cal G}_{{p'}}} = \emptyset, \; \forall p \ne p'. 
\end{equation}
For the sake of notational simplicity,  we use $( \hat{\mathscr{X}}_p, {\hat{\mathscr F}_p}, {\hat{{ \mathscr P}}}_p^{\boldsymbol{\theta }} )$ to denote the statistical model for any ${\bf x}_j$ with $j \in {\mathcal G}_p$.

In general, it is possible that some different sensors employ an identical superquanizer\footnote{Here order is important, thus an identical superquantizer uses the same vector quantizers in the same order.} to convert its observation vector to digital data. 
Each ${\mathcal G}_p$ can be further divided into $M_p$ disjoint nonempty subgroups $\{{\mathcal G}_p^{(m)}\}_{m=1}^{M_p}$ of sensors that use different superquantizers
\begin{equation}
{{\cal G}_p} = \mathop  \cup \limits_{m = 1}^{M_p} {{\cal G}_p^{(m)}}, \text{ and } {{\cal G}_p^{(m)}}\cap {{\cal G}_p^{(m')}} = \emptyset, \; \forall m \ne m'.
\end{equation}
For simplicity, we use 
\begin{equation} \label{hat_Gamma_pm}
\hat \Gamma _p^{(m)} \buildrel \Delta \over = {\left[ {\hat \gamma _{p1}^{(m)},\hat \gamma _{p2}^{(m)},...,\hat \gamma _{p{L_p^{(m)}}}^{(m)}} \right]^T}
\end{equation} 
to denote the superquantizer employed by the sensors in ${{\cal G}_p^{(m)}}$, where $L_p^{(m)}$ is the number of vector quantizers in ${\hat \Gamma }_p^{(m)}$. Moreover, 
we use ${\hat R}_{pl}^{(m)}$ and $\{{\hat{I}}_{mpl}^{(r)}\}_{r=1}^{{\hat R}_{pl}^{(m)}}$ to respecitvely denote the number of quantization levels of ${{\hat \gamma }_{p{l}}^{(m)}}$ and the quantization regions of ${{\hat \gamma }_{p{l}}^{(m)}}$ for each $m$, $p$ and  $l$.

Thus, under \emph{Assumption \ref{Assumption_identical_model}}, if $j$ and $j'$ are contained in ${\mathcal G}_p^{(m)}$ for some $p$ and $m$, then for any $\boldsymbol{\theta}$ and any outcome ${\bf s}$ of the superquantizer $\hat \Gamma _p^{(m)}$,
\begin{equation} \label{qj_qj}
q_j^{({\bf{s}})}\left( {\boldsymbol{\theta }} \right) = q_{{j'}}^{({\bf{s}})}\left( {\boldsymbol{\theta }} \right) = {\hat{{ \mathscr P}}}_p^{\boldsymbol{\theta }}\left( {\hat \Gamma _p^{(m)}\left( {{{\bf{x}}_j}} \right) = {\bf s}} \right).
\end{equation}

By employing (\ref{qj_qj}) and similar arguments as those in \emph{Section \ref{Section_Information_Regularity_Estimation_Capacity}} and \emph{Section \ref{Section_Nonidentifiability_Estimation_Capacity}}, the following theorem can be obtained.
\begin{theorem} \label{Theorem_identical_model}
	Under \emph{Assumption \ref{A1_Open_Set_Domain}}, \emph{Assumption \ref{A2_differentiable}} and \emph{Assumption \ref{Assumption_identical_model}}, for any given $\boldsymbol{\theta}$, any quantization regions $\{{\hat{I}}_{mpl}^{(r)}\}$ and any statistical models $\{ ( \hat{\mathscr{X}}_p, {\hat{\mathscr F}_p}, {\hat{{ \mathscr P}}}_p^{\boldsymbol{\theta }} ) \}$, if the dimension $D_{\boldsymbol{\theta}}$ of the vector parameter $\boldsymbol{\theta}$ is greater than $\sum\nolimits_{p = 1}^P {\sum\nolimits_{m = 1}^{{M_p}} {( {\prod\nolimits_{l = 1}^{L_p^{(m)}} {\hat R_{pl}^{(m)}}  - 1} )} } $, i.e., 
	\begin{align} \notag
	D_{\boldsymbol{\theta}} &  > {\lambda _{{\text{ISM}}}}\left( {\{ {{\cal G}_p^{(m)}} \},\{ {\hat R_{pl}^{(m)}} \}} \right) \\ \label{Define_lambda_identical_model}
	& \buildrel \Delta \over = \sum\limits_{p = 1}^P {\sum\limits_{m = 1}^{{M_p}} {\left( {\prod\limits_{l = 1}^{L_p^{(m)}} {\hat R_{pl}^{(m)}}  - 1} \right)} }, 
	\end{align}
	then the FIM for estimating ${\boldsymbol{\theta}}$ is singular.
	Furthermore, under \emph{Assumption \ref{A1_Open_Set_Domain}}, \emph{Assumption \ref{A3_continuity}} and \emph{Assumption \ref{Assumption_identical_model}}, for any given $\{{\hat{I}}_{mpl}^{(r)}\}$ and $\{ ( \hat{\mathscr{X}}_p, {\hat{\mathscr F}_p}, {\hat{{ \mathscr P}}}_p^{\boldsymbol{\theta }} ) \}$, if (\ref{Define_lambda_identical_model}) holds, 
	then the vector parameter space $\boldsymbol{\Theta}$ is not identifiable. Moreover, for any open subset ${\cal U} \subset {\boldsymbol{\Theta}}$ in ${\mathbbm{R}}^{D_{\boldsymbol{\theta}}}$, there are infinitely many vector parameter points in ${\cal U}$ which are not identifiable.
\end{theorem}

The proof of \emph{Theorem \ref{Theorem_identical_model}} is omitted, since it is similar to the proofs of \emph{Theorem \ref{Theorem_Singularity_vector_quantization}} and \emph{Theorem \ref{Theorem_nonidentifiable_vector_parameter_space}} after 
properly accounting for the impact of \emph{Assumption \ref{Assumption_identical_model}}, which effectively reduces the IDQD. At least in terms of 
the FIM singularity, the fact that the IDQD is reduced seems reasonable since the statistically identical models assumed in \emph{Assumption \ref{Assumption_identical_model}} leads to identical terms in the sum in (\ref{FIM_vector_quantization}) which leads to a smaller dimension 
of ${\boldsymbol{\theta}}$ at which the FIM must be singular.  The impact of \emph{Assumption \ref{Assumption_identical_model}} on identifiability can be similarly justified. 

By the definition of ${\lambda _{{\text{ISM}}}}( {\{ {{\cal G}_p^{(m)}} \},\{ {\hat R_{pl}^{(m)}} \}} )$ in (\ref{Define_lambda_identical_model}), we know that the critical quantity ${\lambda _{{\text{ISM}}}}( {\{ {{\cal G}_p^{(m)}} \},\{ {\hat R_{pl}^{(m)}} \}} )$  does not depend on $\{{\hat{I}}_{mpl}^{(r)}\}$ and $\{ ( \hat{\mathscr{X}}_p, {\hat{\mathscr F}_p}, {\hat{{ \mathscr P}}}_p^{\boldsymbol{\theta }} ) \}$, but is only determined by the number of groups $\{{\mathcal G}_p\}$, the number of subgroups $\{{\mathcal G}_p^{(m)}\}$ and the precision of the vector quantizers employed by the system. 

What's more, we can obtain the relationship between the IDQD $\lambda ( {N,\{ {{R_{jl}}} \}} )$ in (\ref{Estimation_Capacity_Vector_Quantization}) and the quantity ${\lambda _{{\text{ISM}}}}( {\{ {{\cal G}_p^{(m)}} \},\{ {\hat R_{pl}^{(m)}} \}} )$ in (\ref{Define_lambda_identical_model}) as
\begin{align} \notag
\lambda \left( {N,\{ {{R_{jl}}} \}} \right) & = \sum\limits_{j = 1}^N {\prod\limits_{l = 1}^{{L_j}} {{R_{jl}}} }  - N\\ \notag
& = \sum\limits_{p = 1}^P {\sum\limits_{m = 1}^{{M_p}} {\sum\limits_{j \in {\cal G}_p^{(m)}} {\left( {\prod\limits_{l = 1}^{L_p^{(m)}} {\hat R_{pl}^{(m)}}  - 1} \right)} } } \\ \notag
& = \sum\limits_{p = 1}^P {\sum\limits_{m = 1}^{{M_p}} {\left| {{\cal G}_p^{(m)}} \right|\left( {\prod\limits_{l = 1}^{L_p^{(m)}} {\hat R_{pl}^{(m)}}  - 1} \right)} } \\ \label{temp_inequality_identical_model}
& \ge \sum\limits_{p = 1}^P {\sum\limits_{m = 1}^{{M_p}} {\left( {\prod\limits_{l = 1}^{L_p^{(m)}} {\hat R_{pl}^{(m)}}  - 1} \right)} } \\ \label{inequality_identical_model}
& = {\lambda _{{\rm{ISM}}}}\left( {\{ {{\cal G}_p^{(m)}} \},\{ {\hat R_{pl}^{(m)}} \}} \right)
\end{align} 
where the inequality in (\ref{temp_inequality_identical_model}) is a consequence of the fact that $| {{\cal G}_p^{(m)}} | \ge 1$ for all $p$ and all $m$. 
Thus from (\ref{inequality_identical_model}), rather than utilizing the IDQD $\lambda ( {N,\left\{ {{R_{jl}}} \right\}} )$ in (\ref{Estimation_Capacity_Vector_Quantization}), it would be better to employ the critical quantity ${\lambda _{{\text{ISM}}}}( {\{ {{\cal G}_p^{(m)}} \},\{ {\hat R_{pl}^{(m)}} \}} )$ in (\ref{Define_lambda_identical_model}) to specify the fundamental limitation of the quantized estimation system under \emph{Assumption \ref{Assumption_identical_model}} given the conditions of \emph{Theorem \ref{Theorem_identical_model}} apply. 
To this end,  the critical quantity ${\lambda _{{\text{ISM}}}}( {\{ {{\cal G}_p^{(m)}} \},\{ {\hat R_{pl}^{(m)}} \}} )$ in (\ref{Define_lambda_identical_model}) is referred to as the \emph{refined inestimable dimension for quantized data} (rIDQD) for the quantized estimation system under \emph{Assumption \ref{Assumption_identical_model}}.

Additionally, (\ref{temp_inequality_identical_model}) implies that in order to reduce the severity of the fundamental limitation of the quantized estimation system, for any set of sensors whose observation vectors obey the same statistical model, we should employ distinct superquantizers at each of the sensors in this set, so that we can achieve $| {{\cal G}_p^{(m)}} | = 1$ for all $p$ and all $m$. Otherwise, the quantization induced fundamental limitation becomes  more limiting implying the FIM singularity and the nonidentifiability of the vector parameter space for an even smaller vector parameter dimension. 



\subsection{Independent Observation Subvectors}
\label{Section_Assumption_Independent}

In this subsection, we make the following assumption.
\begin{assumption} \label{Assumption_independent}
	All the partitioned observation subvectors $\{{\bf x}_{jl}\}$ are known to be independent. We denote the statistical model of ${\bf x}_{jl}$ by $( {{{\mathscr{X}}_{jl}},{{\mathscr{F}}_{jl}}, {\mathscr{P}}_{jl}^{\boldsymbol{\theta }}} )$ for each $j$ and $l$.
\end{assumption} 

Note that \emph{Assumption \ref{Assumption_independent}} is commonly assumed in recent literature on parameter estimation with quantized data, see \cite{shen2014robust} for example. It is clear that the quantized estimation system under \emph{Assumption \ref{Assumption_independent}} is a special case of the general quantized estimation system described in \emph{Section \ref{Section_System_Model}}.  However, as stated previously we show we can find a smaller IDQD under \emph{Assumption \ref{Assumption_independent}} that better describes the limitations imposed by using quantized data. 

It is clear that under \emph{Assumption \ref{Assumption_independent}}, for any $\boldsymbol{\theta}$ and any outcome ${\bf s}=[s_1,s_2,...,s_{L_j}]^T$ of the superquantizer $ \Gamma _j$, we have
\begin{align} \notag
q_j^{({\bf{s}})}\left( {\boldsymbol{\theta }} \right) & = {\mathscr{P}}_j^{\boldsymbol{\theta }}\left( {{\Gamma _j}\left( {{{\bf{x}}_j}} \right) = {\bf{s}}} \right) \\ \label{q_j_indep}
&  = \prod\limits_{l = 1}^{{L_j}} {{{\mathscr{P}}_{jl}^{\boldsymbol{\theta }}}\left( {{\gamma _{jl}}\left( {{{\bf{x}}_{jl}}} \right) = {s_l}} \right)}. 
\end{align}

Thus, under \emph{Assumption \ref{Assumption_independent}}, we can obtain the following theorem by employing (\ref{q_j_indep}) and similar arguments as those in \emph{Section \ref{Section_Information_Regularity_Estimation_Capacity}} and \emph{Section \ref{Section_Nonidentifiability_Estimation_Capacity}}.
\begin{theorem} \label{Theorem_independent}
	Under \emph{Assumption \ref{A1_Open_Set_Domain}}, \emph{Assumption \ref{A2_differentiable}} and \emph{Assumption \ref{Assumption_independent}}, for any given $\boldsymbol{\theta}$, any quantization regions $\{ {I_{jl}^{(r)}} \}$ and any statistical models $\{ ( {{{\mathscr{X}}_{jl}},{{\mathscr{F}}_{jl}}, {\mathscr{P}}_{jl}^{\boldsymbol{\theta }}} ) \}$, if the dimension $D_{\boldsymbol{\theta}}$ of the vector parameter $\boldsymbol{\theta}$ is greater than $\sum\nolimits_{j = 1}^N {\sum\nolimits_{l = 1}^{{L_j}} {{R_{jl}}} }  - \sum\nolimits_{j = 1}^N {{L_j}} $, i.e., 
	\begin{equation} \label{Define_lambda_independent}
	D_{\boldsymbol{\theta}} > {\lambda _{{\text{Indep}}}}\left( {N,\left\{ {{R_{jl}}} \right\}} \right) \buildrel \Delta \over = \sum\limits_{j = 1}^N {\sum\limits_{l = 1}^{{L_j}} {{R_{jl}}} }  - \sum\limits_{j = 1}^N {{L_j}},
	\end{equation}
	then the FIM for estimating ${\boldsymbol{\theta}}$ is singular.
	Furthermore, under \emph{Assumption \ref{A1_Open_Set_Domain}}, \emph{Assumption \ref{A3_continuity}} and \emph{Assumption \ref{Assumption_independent}}, for any given $\{ {I_{jl}^{(r)}} \}$ and $\{ ( {{{\mathscr{X}}_{jl}},{{\mathscr{F}}_{jl}}, {\mathscr{P}}_{jl}^{\boldsymbol{\theta }}} ) \}$, if (\ref{Define_lambda_independent}) holds, 
	then the vector parameter space $\boldsymbol{\Theta}$ is not identifiable. Moreover, for any open subset ${\cal U} \subset {\boldsymbol{\Theta}}$ in ${\mathbbm{R}}^{D_{\boldsymbol{\theta}}}$, there are infinitely many vector parameter points in ${\cal U}$ which are not identifiable.
\end{theorem}

\emph{Theorem \ref{Theorem_independent}} can be justified as we now explain. Under \emph{Assumption \ref{Assumption_independent}}, all the partitioned observation subvectors $\{{\bf x}_{jl}\}$ are known to be independent. Hence, without any impact, for each $j$, we can view the $j$-th sensor as $L_j$ effective ``sensors'' where the observation vector of the $l$-th effective ``sensor'' is ${\bf x}_{jl}$ and the $l$-th effective ``sensor'' employs the vector quantizer $\gamma_{jl}$ to quantize its observation vector ${\bf x}_{jl}$. As a result, under \emph{Assumption \ref{Assumption_independent}}, the original $N$-sensor system where the $j$-th sensor employs the superquantizer $\Gamma_j$ for each $j$, is equivalent to a $(\sum\nolimits_{j = 1}^N {{L_j}} )$-sensor system where each sensor just employs a vector quantizer. We use a pair of indices $\{(j,l)\}_{j=1,l=1}^{N,L_j}$ to index the sensors in the $(\sum\nolimits_{j = 1}^N {{L_j}} )$-sensor system, and the number of quantization levels of the vector quantizer employed at the $(j,l)$-th sensor is $R_{jl}$. Therefore, by the equivalence between these two sensor systems and by replacing $N$ by $(\sum\nolimits_{j = 1}^N {{L_j}} )$ and replacing $\sum\nolimits_{j = 1}^N {\prod\nolimits_{l = 1}^{L_j} {{R_{jl}}} } $ by $\sum\nolimits_{j = 1}^N {\sum\nolimits_{l = 1}^{{L_j}} {{R_{jl}}} } $ (the new sum over all sensors of the number of quantization levels at each sensor) in the formula in (\ref{Estimation_Capacity_Vector_Quantization}), it follows that 
\begin{equation} \label{explain_lambda_indep}
{\lambda _{{\text{Indep}}}}\left( {N,\left\{ {{R_{jl}}} \right\}} \right) = \sum\limits_{j = 1}^N {\sum\limits_{l = 1}^{{L_j}} {{R_{jl}}} }  - \sum\limits_{j = 1}^N {{L_j}},
\end{equation}
which justifies \emph{Theorem \ref{Theorem_independent}}.

Noting that $R_{jl} \ge 1$ for all $j$ and $l$, and by employing the fact that for any positive integer $L_j$, if $a_i \ge 0$ for all $i=1,2,..,L_j$, then $\prod\nolimits_{i = 1}^{L_j} {\left( {1 + {a_i}} \right)}  \ge 1 + \sum\nolimits_{i = 1}^{L_j} {{a_i}}$, we can obtain
\begin{align} \notag
\lambda \left( {N,\left\{ {{R_{jl}}} \right\}} \right) & = \sum\limits_{j = 1}^N {\prod\limits_{l = 1}^{{L_j}} {{R_{jl}}} }  - N\\  \notag
& = \sum\limits_{j = 1}^N {\left\{ {\prod\limits_{l = 1}^{{L_j}} {\left[ {1 + \left( {{R_{jl}} - 1} \right)} \right]}  - 1} \right\}} \\ \notag
& \ge \sum\limits_{j = 1}^N {\left\{ {1 + \sum\limits_{l = 1}^{{L_j}} {\left( {{R_{jl}} - 1} \right)}  - 1} \right\}} \\ \label{inequality_independent}
&= {\lambda _{{\text{Indep}}}}\left( {N,\left\{ {{R_{jl}}} \right\}} \right).
\end{align}
Thus, it is seen from (\ref{inequality_independent}) that the critical quantity ${\lambda _{{\text{Indep}}}}\left( {N,\left\{ {{R_{jl}}} \right\}} \right)$ in (\ref{Define_lambda_independent}) 
allows us to guarantee the FIM is singular and the vector parameter space is nonidentifiable at a smaller dimension of $\boldsymbol{\theta}$. We refer to ${\lambda _{{\text{Indep}}}}\left( {N,\left\{ {{R_{jl}}} \right\}} \right)$ in (\ref{Define_lambda_independent}) as the rIDQD for the quantized estimation system under \emph{Assumption \ref{Assumption_independent}}. 

What's more, it is seen from (\ref{Define_lambda_independent}) that the rIDQD ${\lambda _{{\text{Indep}}}}\left( {N,\left\{ {{R_{jl}}} \right\}} \right)$ under \emph{Assumption \ref{Assumption_independent}} is precisely the number of quantization levels employed by the effective sensor system minus the number of effective sensors. Hence, the rIDQD ${\lambda _{{\text{Indep}}}}\left( {N,\left\{ {{R_{jl}}} \right\}} \right)$ under \emph{Assumption \ref{Assumption_independent}} does not depend on $\{I_{jl}^{(r)}\}$ and $\{ ( {{{\mathscr{X}}_{jl}},{{\mathscr{F}}_{jl}}, {\mathscr{P}}_{jl}^{\boldsymbol{\theta }}} ) \}$.

\subsection{Independent Observation Subvectors, Identical Sensor Observation Statistical Models and Identical Vector Quantizers}

In this subsection, we consider the following assumption which is  stronger than \emph{Assumption \ref{Assumption_independent}}.
\begin{assumption} \label{Assumption_indep_identical}
	All the partitioned observation subvectors $\{{\bf x}_{jl}\}$ are known to be independent, and moreover, some partitioned observation subvectors ${\bf x}_{jl}$ are known to obey the same statistical model such that $W<\sum\nolimits_{j = 1}^N {{L_j}} $ distinct statistical models of the partitioned observation subvectors $\{{\bf x}_{jl}\}$ exist.
\end{assumption}


Collect all the partitioned observation subvector indices that employ the $w$-th statistical model in the group ${{\mathcal A}_w}$.
%
For notational simplicity,  we use $( {{\tilde{\mathscr{X}}_{w}},{\tilde{\mathscr{F}}_{w}}, \tilde{\mathscr{P}}_{w}^{\boldsymbol{\theta }}} ) $  to denote the statistical model for any ${\bf x}_{jl}$ with its index contained in ${\mathcal A}_w$.

Moreover, each group ${\mathcal A}_w$ can be divided into $T_w$ disjoint nonempty subgroups $\{{\mathcal A}_w^{(t)}\}_{t=1}^{T_w}$ of partitioned observation subvector indices that employ different vector quantizers such that
\begin{equation}
{{\cal A}_w} = \mathop  \cup \limits_{t = 1}^{T_w} {{\cal A}_w^{(t)}}, \text{ and } {{\cal A}_w^{(t)}}\cap {{\cal A}_w^{(t')}} = \emptyset, \; \forall t \ne t', 
\end{equation}
In other words, if the indices of ${\bf x}_{j_1 l_1}$ and  ${\bf x}_{j_2 l_2}$ are contained in some ${\mathcal A}_w^{(t)}$, then $\gamma_{j_1 l_1} = \gamma_{j_2 l_2}$. For simplicity, we use $\tilde{\gamma}_w^{(t)}$ to denote the vector quantizer employed for the partitioned observation subvectors whose indices are contained in ${{\cal A}_w^{(t)}}$, and use ${\tilde R}_{w}^{(t)}$ and $\{{\tilde{I}}_{wt}^{(r)}\}_{r=1}^{{\tilde R}_{w}^{(t)}}$ to respectively denote the number of quantization levels of $\tilde{\gamma}_w^{(t)}$ and the quantization regions of $\tilde{\gamma}_w^{(t)}$ for each $w$ and $t$.

Under \emph{Assumption \ref{Assumption_indep_identical}}, we have the following theorem with regard to the fundamental limitation of the quantized estimation system.
\begin{theorem} \label{Theorem_indep_identical}
	Under \emph{Assumption \ref{A1_Open_Set_Domain}}, \emph{Assumption \ref{A2_differentiable}} and \emph{Assumption \ref{Assumption_indep_identical}}, for any given $\boldsymbol{\theta}$, any quantization regions $\{{\tilde{I}}_{wt}^{(r)}\}$ and any statistical models $\{ ( {{\tilde{\mathscr{X}}_{w}},{\tilde{\mathscr{F}}_{w}}, \tilde{\mathscr{P}}_{w}^{\boldsymbol{\theta }}} ) \}$, if the dimension $D_{\boldsymbol{\theta}}$ of the vector parameter $\boldsymbol{\theta}$ is greater than $\sum\nolimits_{w = 1}^W {\sum\nolimits_{t = 1}^{{T_w}} {\left( {\tilde R_w^{(t)} - 1} \right)} } $, i.e., 
	\begin{equation}  \label{Define_lambda_indep_identical}
	D_{\boldsymbol{\theta}}   >  \lambda _{{\text{Indep}}}^{{\text{ISM}}}\left( {\{ {{\cal A}_w^{(t)}} \},\{ {\tilde R_w^{(t)}} \}} \right) \buildrel \Delta \over =  \sum\limits_{w = 1}^W {\sum\limits_{t = 1}^{{T_w}} {\left( {\tilde R_w^{(t)} - 1} \right)} },
	\end{equation}
	then the FIM for estimating ${\boldsymbol{\theta}}$ is singular.
	Furthermore, under \emph{Assumption \ref{A1_Open_Set_Domain}}, \emph{Assumption \ref{A3_continuity}} and \emph{Assumption \ref{Assumption_indep_identical}}, for any given $\{{\tilde{I}}_{wt}^{(r)}\}$ and $\{ ( {{\tilde{\mathscr{X}}_{w}},{\tilde{\mathscr{F}}_{w}}, \tilde{\mathscr{P}}_{w}^{\boldsymbol{\theta }}} ) \}$, if (\ref{Define_lambda_indep_identical}) holds, 
	then the vector parameter space $\boldsymbol{\Theta}$ is not identifiable. Moreover, for any open subset ${\cal U} \subset {\boldsymbol{\Theta}}$ in ${\mathbbm{R}}^{D_{\boldsymbol{\theta}}}$, there are infinitely many vector parameter points in ${\cal U}$ which are not identifiable.
\end{theorem}

Since \emph{Assumption \ref{Assumption_indep_identical}} combines \emph{Assumption \ref{Assumption_identical_model}} and \emph{Assumption \ref{Assumption_independent}}, the proof of \emph{Theorem \ref{Theorem_indep_identical}} involves a combination of the proofs of \emph{Theorem \ref{Theorem_identical_model}} and \emph{Theorem \ref{Theorem_independent}}.


Note that
\begin{align} \notag
{\lambda _{{\text{Indep}}}}\left( {N,\left\{ {{R_{jl}}} \right\}} \right) & = \sum\limits_{j = 1}^N {\sum\limits_{l = 1}^{{L_j}} {\left( {{R_{jl}} - 1} \right)} } \\ \notag
& = \sum\limits_{w = 1}^W {\sum\limits_{t = 1}^{{T_w}} {\left| {{\cal A}_w^{(t)}} \right|\left( {\tilde R_w^{(t)} - 1} \right)} } \\ \label{temp_inequality_indep_identical}
& \ge \sum\limits_{w = 1}^W {\sum\limits_{t = 1}^{{T_w}} {\left( {\tilde R_w^{(t)} - 1} \right)} } \\  \label{inequality_indep_identical}
&= \lambda _{{\text{Indep}}}^{{\text{ISM}}}\left( {\{ {{\cal A}_w^{(t)}} \},\{ {\tilde R_w^{(t)}} \}} \right),
\end{align}
where (\ref{temp_inequality_indep_identical}) is based on the fact that $| {{\cal A}_w^{(t)}} | \ge 1$ for all $w$ and all $t$. Therefore, under \emph{Assumption \ref{Assumption_indep_identical}} which is stronger than \emph{Assumption \ref{Assumption_independent}}, the sufficient condition in (\ref{Define_lambda_independent}) is even less restrictive than the sufficient condition in (\ref{Define_lambda_indep_identical}) which considers scenarios under \emph{Assumption \ref{Assumption_independent}}.
We call the quantity $\lambda _{{\text{Indep}}}^{{\text{ISM}}}( {\{ {{\cal A}_w^{(t)}} \},\{ {\tilde R_w^{(t)}} \}} )$ in (\ref{Define_lambda_indep_identical}) the rIDQD for the quantized estimation system under \emph{Assumption \ref{Assumption_indep_identical}}.

It should be noted that under some other assumptions, we can also obtain the corresponding rIDQD by employing  similar arguments to those just presented. For the sake of brevity, we omit the detailed discussion.

\section{Conclusion} 
\label{Section_Conclusion}

In this paper, we investigate the impact of quantization on the estimation capabilities with respect to
the information-regularity condition and the identifiability condition. A critical quantity, called IDQD,
is introduced, which describes a fundamental limitation of using quantized data. To be specific, under
the condition that the dimension of the desired vector parameter is larger than the IDQD, the FIM for
estimating the desired vector parameter is singular for any value of the desired vector parameter, any
quantization regions, and any statistical models of the observations. Furthermore, it is shown that under
the same condition, the vector parameter space is not identifiable, and moreover, there are infinitely
many nonidentifiable vector parameter points in the vector parameter space. It is worth mentioning that
there is no general equivalence between the quantization induced FIM singularity and the quantization
induced nonidentifiability of the vector parameter space. Further, in the quantization induced nonidentifiable vector parameter space, every vector parameter point is nonidentifiable in some cases, while in some other cases, there exist some identifiable vector parameter points. Thus the quantization induced FIM singularity does not necessarily determine the identifiability of the vector parameter point although it does determine
the identifiability of the vector parameter space. Moreover, the cardinality of a set of observationally
equivalent points in the quantization induced nonidentifiable vector parameter space can be as small
as 1 and can also be as large as uncountably infinite. In addition, some commonly assumed specific
assumptions on the statistical models of the observations are considered in this paper. It is shown that
under these assumptions, a refined IDQD becomes smaller than the standard  IDQD, implying the FIM singularity and the nonidentifiability of the vector parameter space can be guaranteed for an even smaller
vector parameter dimension.

\appendices

\section{Proof of Theorem \ref{Theorem_Singularity_vector_quantization}} \label{proof_Theorem_Singularity_vector_quantization}
By employing (\ref{FIM_vector_quantization}), the rank of ${\bf{J}}\left(  {\boldsymbol{\theta}}   \right) $ is upper bounded by
\begin{align} \notag
& {\rm{rank}}\left( {{\bf{J}}\left( {\boldsymbol{\theta }} \right)} \right)\\ \notag
& = {\rm{rank}}\left( {\sum\limits_{j = 1}^N {\sum\limits_{{\bf{s}} \in {{\cal S}_j}} {\frac{1}{{q_j^{({\bf{s}})}\left( {\boldsymbol{\theta }} \right)}}\frac{{\partial q_j^{({\bf{s}})}\left( {\boldsymbol{\theta }} \right)}}{{\partial {\boldsymbol{\theta }}}}{{\left[ {\frac{{\partial q_j^{({\bf{s}})}\left( {\boldsymbol{\theta }} \right)}}{{\partial {\boldsymbol{\theta }}}}} \right]}^T}} } } \right)\\ \label{proof_Theorem_Singularity_vector_quantization_temp1}
& \le \sum\limits_{j = 1}^N {{\rm{rank}}\left( {\sum\limits_{{\bf{s}} \in {{\cal S}_j}} {\frac{{\partial q_j^{({\bf{s}})}\left( {\boldsymbol{\theta }} \right)}}{{\partial {\boldsymbol{\theta }}}}{{\left[ {\frac{{\partial q_j^{({\bf{s}})}\left( {\boldsymbol{\theta }} \right)}}{{\partial {\boldsymbol{\theta }}}}} \right]}^T}} } \right)}.
\end{align}
Noticing that 
\begin{equation}
\sum\limits_{{\bf{s}} \in {{\cal S}_j}} {q_j^{({\bf{s}})}\left( {\boldsymbol{\theta }} \right)}  = 1, \; \forall j,
\end{equation}
we can obtain that 
\begin{equation}
\sum\limits_{{\bf{s}} \in {{\cal S}_j}} {\frac{{\partial q_j^{({\bf{s}})}\left( {\boldsymbol{\theta }} \right)}}{{\partial {\boldsymbol{\theta }}}}}  = {\bf 0},  \; \forall j,
\end{equation}
and therefore, 
\begin{align} \notag
& {\rm{rank}}\left( {\sum\limits_{{\bf{s}} \in {{\cal S}_j}} {\frac{{\partial q_j^{({\bf{s}})}\left( {\boldsymbol{\theta }} \right)}}{{\partial {\boldsymbol{\theta }}}}{{\left[ {\frac{{\partial q_j^{({\bf{s}})}\left( {\boldsymbol{\theta }} \right)}}{{\partial {\boldsymbol{\theta }}}}} \right]}^T}} } \right) \\ \notag
& \le \left| {{{\cal S}_j}} \right| - 1 \\ \label{proof_Theorem_Singularity_vector_quantization_temp2}
& = \prod\limits_{l = 1}^{{L_j}} {{R_{jl}}}  - 1, \; \forall j,
\end{align}
where (\ref{proof_Theorem_Singularity_vector_quantization_temp2}) follows from (\ref{size_Sp}).

By employing (\ref{proof_Theorem_Singularity_vector_quantization_temp1}) and (\ref{proof_Theorem_Singularity_vector_quantization_temp2}), we can bound the rank of ${\bf{J}}\left(  {\boldsymbol{\theta}}   \right) $ above by
\begin{align} \notag
{\rm{rank}}\left( {{\bf{J}}\left( {\boldsymbol{\theta }} \right)} \right) & \le \sum\limits_{j = 1}^N { {\left( {\prod\limits_{l = 1}^{{L_j}} {{R_{jl}}}  - 1} \right)} } \\ \label{proof_Theorem_Singularity_vector_quantization_temp3}
& = \sum\limits_{j = 1}^N {\prod\limits_{l = 1}^{{L_j}} {{R_{jl}}} }  - N.
\end{align}

Thus, noting that the size of ${\bf{J}}\left(  {\boldsymbol{\theta}}   \right) $ is  $D_{\boldsymbol{\theta}}$-by-$D_{\boldsymbol{\theta}}$, if $D_{\boldsymbol{\theta}} > \lambda \left( {N,\left\{ {{R_{jl}}} \right\}} \right)  = \sum_{j = 1}^N {\prod_{l = 1}^{{L_j}} {{R_{jl}}} }  - N$, ${\bf{J}}\left(  {\boldsymbol{\theta}}   \right) $  is singular for any given $\boldsymbol{\theta}$, $\{ {I_{jl}^{(r)}} \}$ and $\{( {\mathscr{X}}_j, {{\mathscr F}_j}, {{{ \mathscr P}}}_j^{\boldsymbol{\theta }} )\}$.

\section{Proof of Lemma \ref{Lemma_iff_identifiability_condition}} \label{proof_Lemma_iff_identifiability_condition}

First, consider a mapping ${{\boldsymbol{\bar \Psi }}}$
\begin{equation} \label{map_bar_Psi}
\begin{aligned}  
{{\boldsymbol{\bar \Psi }}} :   {\boldsymbol{\Theta }} & \longrightarrow {{\mathbbm{R}}^{{\sum\limits_{j = 1}^N {\prod\limits_{l = 1}^{{L_j}} {{R_{jl}}} } }}}\\ 
{\boldsymbol{\theta}} &  \longmapsto {{\boldsymbol{\bar \Psi }}}\left( {\boldsymbol{\theta }} \right),
\end{aligned}
\end{equation}
where the $(\sum\nolimits_{j = 1}^N {\prod\nolimits_{l = 1}^{{L_j}} {{R_{jl}}} } )$-dimensional vector  ${{\boldsymbol{\bar \Psi }}}\left( {\boldsymbol{\theta }} \right)$ is defined as
\begin{equation} \label{bar_Psi}
{\boldsymbol{\bar \Psi }}\left( {\boldsymbol{\theta }} \right) \buildrel \Delta \over = {\left[ {{{{\boldsymbol{\bar \psi }}}_1}{{\left( {\boldsymbol{\theta }} \right)}^T},{{{\boldsymbol{\bar \psi }}}_2}{{\left( {\boldsymbol{\theta }} \right)}^T},...,{{{\boldsymbol{\bar \psi }}}_N}{{\left( {\boldsymbol{\theta }} \right)}^T}} \right]^T},
\end{equation}
and for each $j$, the $|{{\cal S}_j}|$-dimensional vector ${{{\boldsymbol{\bar \psi }}}_j}\left( {\boldsymbol{\theta }} \right)$ is defined as
\begin{equation} \label{bar_psi}
{{{\boldsymbol{\bar \psi }}}_j}\left( {\boldsymbol{\theta }} \right) \buildrel \Delta \over = {\left[ {q_j^{({\bf{s}}_1^{(j)})}\left( {\boldsymbol{\theta }} \right),q_j^{({\bf{s}}_2^{(j)})}\left( {\boldsymbol{\theta }} \right),...,q_j^{({\bf{s}}_{|{{\cal S}_j}|}^{(j)})}\left( {\boldsymbol{\theta }} \right)} \right]^T},
\end{equation}
$q_{j}^{({\bf s})}\left( {\boldsymbol{\theta }} \right)$ is defined in (\ref{q_j_s}), and ${\bf s}_i^{(j)}$ is defined in (\ref{S_p}) for all $i=1,2,...,|{{\cal S}_j}|$.

We first show that the mapping ${\varphi}_{\bf u}$ in (\ref{phi_u}) is not injective if and only if the mapping ${\boldsymbol{\bar \Psi }}$ in (\ref{map_bar_Psi}) is not injective, and hence, the injectivity of the mapping ${\boldsymbol{\bar \Psi }}$ in (\ref{map_bar_Psi}) is the same as that of the mapping ${\varphi}_{\bf u}$ in (\ref{phi_u}).

Suppose the mapping ${\boldsymbol{\bar \Psi }}$ in (\ref{map_bar_Psi}) is not injective. Then, there exist two distinct ${\boldsymbol{\theta}}_1,\; {\boldsymbol{\theta}}_1 \in {\boldsymbol{\Theta}}$ such that ${\boldsymbol{\bar \Psi }}\left( {\boldsymbol{\theta }}_1 \right) = {\boldsymbol{\bar \Psi }}\left( {\boldsymbol{\theta }}_2 \right) $. 

Noting that
\begin{equation} \label{Pr_u}
\Pr \left( {{\bf{u}}\left| {\boldsymbol{\theta }} \right.} \right) = \prod\limits_{j = 1}^N {\prod\limits_{{\bf{s}} \in {{\cal S}_j}} {{{\left[ {q_j^{({\bf{s}})}\left( {\boldsymbol{\theta }} \right)} \right]}^{{\mathbbm{1}}\left\{ {{{\bf{u}}_j} = {\bf{s}}} \right\}}}} } ,
\end{equation}
it is clear that $\Pr \left( {{\bf{u}}\left| {\boldsymbol{\theta }}_1 \right.} \right) = \Pr \left( {{\bf{u}}\left| {{{\boldsymbol{\theta }}_2}} \right.} \right)$ for all $\bf{u} \in {\mathcal{A}}$, if there exist two distinct ${\boldsymbol{\theta}}_1,\; {\boldsymbol{\theta}}_1 \in {\boldsymbol{\Theta}}$ such that ${\boldsymbol{\bar \Psi }}\left( {\boldsymbol{\theta }}_1 \right) = {\boldsymbol{\bar \Psi }}\left( {\boldsymbol{\theta }}_2 \right) $. Thus, the mapping ${\varphi}_{\bf u}$ in (\ref{phi_u}) is not injective.

On the other hand, suppose the mapping ${\varphi}_{\bf u}$ in (\ref{phi_u}) is not injective. Then, there exist two distinct ${\boldsymbol{\theta}}_1,\; {\boldsymbol{\theta}}_1 \in {\boldsymbol{\Theta}}$ such that 
\begin{equation} \label{proof_lemma1_temp}
\Pr \left( {{\bf{u}}\left| {\boldsymbol{\theta }}_1 \right.} \right) = \Pr \left( {{\bf{u}}\left| {{{\boldsymbol{\theta }}_2}} \right.} \right), \; \forall \bf{u} \in {\mathcal{A}},
\end{equation}
where $\Pr \left( {{\bf{u}}\left| {\boldsymbol{\theta }} \right.} \right)$ is defined in (\ref{Pr_u}).

Note that for each $i$, we have
\begin{equation} \label{proof_lemma1_temp3}
\sum\limits_{{\bf{s}} \in {{\cal S}_i}} {q_i^{({\bf{s}})}\left( {\boldsymbol{\theta }}_1 \right)}  = 1.
\end{equation}
Hence, for each $i$, there exists some ${{\bf{h}}_i} \in {\mathcal S}_i$ such that
\begin{equation} \label{proof_lemma1_temp1}
q_i^{({\bf h}_i)}\left( {\boldsymbol{\theta }}_1 \right) \ne 0.
\end{equation} 

For any given $j$ and any given ${\bf s} \in {\mathcal{S}}_j$, we are going to show that $q_j^{({\bf{s}})}( {{{\boldsymbol{\theta }}_1}} ) = q_j^{({\bf{s}})}( {{{\boldsymbol{\theta }}_2}} )$.

Consider a realization of $\bf u$ that 
\begin{equation}
{\bf{u}} = {\left[ {{{\bf{h}}_1^T},{{\bf{h}}_2^T},...,{{\bf{h}}_{j - 1}^T},{\bf{s}}^T,{{\bf{h}}_{j + 1}^T},...,{{\bf{h}}_N^T}} \right]^T}.
\end{equation}
By employing (\ref{Pr_u}) and (\ref{proof_lemma1_temp}), we have
\begin{equation} \label{proof_lemma1_temp4}
q_j^{({\bf{s}})}\left( {{{\boldsymbol{\theta }}_1}} \right)\prod\limits_{i \ne j} {q_i^{({{\bf{h}}_i})}\left( {{{\boldsymbol{\theta }}_1}} \right)}  = q_j^{({\bf{s}})}\left( {{{\boldsymbol{\theta }}_2}} \right)\prod\limits_{i \ne j} {q_i^{({{\bf{h}}_i})}\left( {{{\boldsymbol{\theta }}_2}} \right)},
\end{equation}
By (\ref{proof_lemma1_temp1}), we know that $\prod\nolimits_{i \ne j} {q_i^{({{\bf{h}}_i})}\left( {{{\boldsymbol{\theta }}_1}} \right)}  \ne 0$, and therefore, from (\ref{proof_lemma1_temp4}), we can obtain
\begin{equation} \label{proof_lemma1_temp2}
q_j^{({\bf{s}})}\left( {{{\boldsymbol{\theta }}_1}} \right) = q_j^{({\bf{s}})}\left( {{{\boldsymbol{\theta }}_2}} \right)\frac{{\prod\limits_{i \ne j} {q_i^{({{\bf{h}}_i})}\left( {{{\boldsymbol{\theta }}_2}} \right)} }}{{\prod\limits_{i \ne j} {q_i^{({{\bf{h}}_i})}\left( {{{\boldsymbol{\theta }}_1}} \right)} }}.
\end{equation}
Furthermore, by  noting that $\sum\nolimits_{{\bf{s}} \in {{\cal S}_j}} {q_j^{({\bf{s}})}\left( {{{\boldsymbol{\theta }}_2}} \right)}  = 1$ and  employing (\ref{proof_lemma1_temp3}) and (\ref{proof_lemma1_temp2}), we can obtain
\begin{align} \notag
\frac{{\prod\limits_{i \ne j} {q_i^{({{\bf{h}}_i})}\left( {{{\boldsymbol{\theta }}_2}} \right)} }}{{\prod\limits_{i \ne j} {q_i^{({{\bf{h}}_i})}\left( {{{\boldsymbol{\theta }}_1}} \right)} }} & = \frac{{\prod\limits_{i \ne j} {q_i^{({{\bf{h}}_i})}\left( {{{\boldsymbol{\theta }}_2}} \right)} }}{{\prod\limits_{i \ne j} {q_i^{({{\bf{h}}_i})}\left( {{{\boldsymbol{\theta }}_1}} \right)} }}\sum\limits_{{\bf{s}} \in {{\cal S}_j}} {q_j^{({\bf{s}})}\left( {{{\boldsymbol{\theta }}_2}} \right)} \\ \notag
& = \sum\limits_{{\bf{s}} \in {{\cal S}_j}} {q_j^{({\bf{s}})}\left( {{{\boldsymbol{\theta }}_1}} \right)} \\
& = 1,
\end{align}
which implies
\begin{equation}
q_j^{({\bf{s}})}\left( {{{\boldsymbol{\theta }}_1}} \right) = q_j^{({\bf{s}})}\left( {{{\boldsymbol{\theta }}_2}} \right).
\end{equation}
Therefore, by the definitions of ${\bf{\bar \Psi }}\left( {\boldsymbol{\theta }} \right)$ and ${{{\boldsymbol{\bar \psi }}}_j}\left( {\boldsymbol{\theta }} \right)$ in (\ref{bar_Psi}) and (\ref{bar_psi}), we know that
\begin{equation}
{\bf{\bar \Psi }}\left( {\boldsymbol{\theta }}_1 \right) = {\bf{\bar \Psi }}\left( {\boldsymbol{\theta }}_2 \right), 
\end{equation}
and hence, the mapping ${\boldsymbol{\bar \Psi }}$ in (\ref{map_bar_Psi}) is not injective. As a result, we know that the mapping ${\varphi}_{\bf u}$ in (\ref{phi_u}) is not injective if and only if the mapping ${\boldsymbol{\bar \Psi }}$ in (\ref{map_bar_Psi}) is not injective, which implies that the injectivity of the mapping ${\boldsymbol{\bar \Psi }}$ in (\ref{map_bar_Psi}) is the same as that of the mapping ${\varphi}_{\bf u}$ in (\ref{phi_u}).

Furthermore, for any given $j$, by the definitions of ${{\boldsymbol{\psi }}_j}\left( {\boldsymbol{\theta }} \right)$ and ${{{\boldsymbol{\bar \psi }}}_j}\left( {\boldsymbol{\theta }} \right)$ in (\ref{psi}) and (\ref{bar_psi}) respectively, and noticing that $\sum\nolimits_{{\bf{s}} \in {{\cal S}_j}} {q_j^{({\bf{s}})}\left( {{{\boldsymbol{\theta }}}} \right)}  = 1$ for all $j$, we can express ${{{\boldsymbol{\bar \psi }}}_j}\left( {\boldsymbol{\theta }} \right)$ as
\begin{equation} \label{proof_lemma1_bar_psi}
{{{\boldsymbol{\bar \psi }}}_j}\left( {\boldsymbol{\theta }} \right) = {\left[ {{{\boldsymbol{\psi }}_j}{{\left( {\boldsymbol{\theta }} \right)}^T},1 - {{\bf{1}}^T}{{\boldsymbol{\psi }}_j}\left( {\boldsymbol{\theta }} \right)} \right]^T}.
\end{equation}

It is clear that if ${{{\boldsymbol{\bar \psi }}}_j}\left( {\boldsymbol{\theta }} \right)$ is not injective, then ${{\boldsymbol{\psi }}_j}\left( {\boldsymbol{\theta }} \right)$ is not injective. On the other hand, if ${{\boldsymbol{\psi }}_j}\left( {\boldsymbol{\theta }} \right)$ is not injective, then there exist two distinct  ${\boldsymbol{\theta}}_1,\; {\boldsymbol{\theta}}_1 \in {\boldsymbol{\Theta}}$ such that ${{\boldsymbol{\psi }}_j}\left( {\boldsymbol{\theta }}_1 \right) = {{\boldsymbol{\psi }}_j}\left( {\boldsymbol{\theta }}_2 \right)$, and hence $1 - {{\bf{1}}^T}{{\boldsymbol{\psi }}_j}\left( {\boldsymbol{\theta }}_1 \right) = 1 - {{\bf{1}}^T}{{\boldsymbol{\psi }}_j}\left( {\boldsymbol{\theta }}_2 \right)$. Consequently, we have ${{{\boldsymbol{\bar \psi }}}_j}\left( {\boldsymbol{\theta }}_1 \right) = {{{\boldsymbol{\bar \psi }}}_j}\left( {\boldsymbol{\theta }}_2 \right)$ by (\ref{proof_lemma1_bar_psi}), which implies that ${{{\boldsymbol{\bar \psi }}}_j}\left( {\boldsymbol{\theta }} \right)$ is not injective. Therefore, the injectivity of ${{{\boldsymbol{\bar \psi }}}_j}\left( {\boldsymbol{\theta }} \right)$ is the same as that of ${{\boldsymbol{\psi }}_j}\left( {\boldsymbol{\theta }} \right)$ for all $j$, which implies that ${{{\boldsymbol{\bar \Psi }}}}\left( {\boldsymbol{\theta }} \right)$ is injective if and only if ${{{\boldsymbol{ \Psi }}}}\left( {\boldsymbol{\theta }} \right)$ is injective. Since we have proven that the mapping ${\varphi}_{\bf u}$ in (\ref{phi_u}) is not injective if and only if the mapping ${\boldsymbol{\bar \Psi }}$ in (\ref{map_bar_Psi}) is not injective, we know that the mapping ${\varphi}_{\bf u}$ in (\ref{phi_u}) is injective if and only if the mapping ${\boldsymbol{ \Psi }}$ in (\ref{map_reduce_dimension_Psi}) is injective.

In order to show that the dimension of the vector  ${\boldsymbol{\Psi }}\left( {\boldsymbol{\theta }} \right)$ in (\ref{Psi}) is strictly smaller than that of ${{{\varphi}_{\bf u}} }\left( {\boldsymbol{\theta }} \right)$ in (\ref{phi_u}) for any given $N$ and $\{R_{jl}\}$, it suffices to show that 
\begin{equation} \label{proof_lemma1_dimension}
{D_{\bf{u}}} = \prod\limits_{j = 1}^N {\prod\limits_{l = 1}^{{L_j}} {{R_{jl}}} }  > \sum\limits_{j = 1}^N {\prod\limits_{l = 1}^{{L_j}} {{R_{jl}}} }  - N,
\end{equation}
for any given $N$ and $\{R_{jl}\}$.

Since $N$ is the number of sensors and $R_{jl}$ denotes the number of quantization levels of the quantizer $\gamma_{jl}$ for each $j$ and $l$, we know that $N \ge 1$ and $R_{jl} \ge 1$ for all $j$ and all $l$. 
Hence, we can obtain that
\begin{equation} \label{proof_lemma1_temp_ge1}
\prod\limits_{l = 1}^{{L_j}} {{R_{jl}}}  \ge 1, \; \forall j.
\end{equation}

Furthermore, notice that if $x \ge 0$ and $y \ge 0$, then we have the following inequality
\begin{equation} 
\left( {1 + x} \right)\left( {1 + y} \right) = 1 + x + y + xy \ge 1 + x + y.
\end{equation}  
Therefore, by induction, we can obtain that if $x_i \ge 0$ for all $i=1,2,...,N$, then
\begin{equation} \label{proof_lemma1_claim2}
\prod\limits_{i = 1}^N {\left( {1 + {x_i}} \right)}  \ge 1 + \sum\limits_{i = 1}^N {{x_i}}.
\end{equation}

By employing (\ref{proof_lemma1_temp_ge1}) and (\ref{proof_lemma1_claim2}), we can obtain
\begin{align} \notag
\prod\limits_{j = 1}^N {\prod\limits_{l = 1}^{{L_j}} {{R_{jl}}} }  & = \prod\limits_{j = 1}^N {\left[ {\left( {\prod\limits_{l = 1}^{{L_j}} {{R_{jl}}}  - 1} \right) + 1} \right]} \\ \notag
& \ge 1 + \sum\limits_{j = 1}^N {\left( {\prod\limits_{l = 1}^{{L_j}} {{R_{jl}}}  - 1} \right)} \\  \label{proof_lemma1_temp_explain}
& > \sum\limits_{j = 1}^N {\prod\limits_{l = 1}^{{L_j}} {{R_{jl}}} }  - N.
\end{align}
This completes the proof.

\section{Proof of Theorem \ref{Theorem_nonidentifiable_vector_parameter_space}} \label{proof_Theorem_nonidentifiable_vector_parameter_space}
Under \emph{Assumption \ref{A1_Open_Set_Domain}}, the interior of ${\boldsymbol{\Theta}}$ is not empty. Thus, there exists a subset ${\mathcal{U}}$ of ${\boldsymbol{\Theta}}$ which is open in ${\mathbbm{R}}^{D_{\boldsymbol{\theta}}}$. 

Define a $D_{\boldsymbol{\theta}}$-dimensional vector ${\bf{\hat \Psi }}\left( {\boldsymbol{\theta }} \right)$
\begin{equation} \label{hat_Psi}
{\bf{\hat \Psi }}\left( {\boldsymbol{\theta }} \right) \buildrel \Delta \over = {\left[ {{\bf{\Psi }}{{\left( {\boldsymbol{\theta }} \right)}^T},{{\bf{0}}^T}} \right]^T},
\end{equation}
where ${\bf{\Psi }}{{\left( {\boldsymbol{\theta }} \right)}}$ is defined in (\ref{Psi}), and the dimension of the all-zero vector in (\ref{hat_Psi}) is ${D_{\boldsymbol{\theta }}} - \lambda \left( {N,\left\{ {{R_{jl}}} \right\}} \right)$.

Under \emph{Assumption \ref{A3_continuity}}, for all $j$ and all ${\bf s}$, ${q_{j}^{({\bf s})}\left( {\boldsymbol{\theta}} \right)}$ is a continuous function with respect to $\boldsymbol{\theta}$. Hence, by (\ref{Psi}), (\ref{psi}) and (\ref{hat_Psi}), we know that the restriction ${\bf{\hat \Psi }}_{\upharpoonright \cal U} $ of the mapping ${\bf{\hat \Psi }}$ to $\cal U$
\begin{equation}\label{map_hat_Psi_U}
\begin{aligned} 
{\bf{\hat \Psi }}_{\upharpoonright \cal U} :  {\cal U} & \longrightarrow {{\mathbbm{R}}^{D_{\boldsymbol{\theta}}}}\\ 
{\boldsymbol{\theta}} &  \longmapsto {\bf{\hat \Psi }}\left( {\boldsymbol{\theta }} \right)
\end{aligned}
\end{equation}
is continuous with respect to $\boldsymbol{\theta}$.

It is clear that ${\mathcal{U}}$ is an open set in ${{\mathbbm{R}}^{D_{\boldsymbol{\theta}}}}$, but by the definition of ${\bf{\hat \Psi }}( {\boldsymbol{\theta }} )$ in (\ref{hat_Psi}), ${\bf{\hat \Psi }}_{\upharpoonright \cal U}\left( {\mathcal{U}} \right)$ is not open in ${{\mathbbm{R}}^{D_{\boldsymbol{\theta}}}}$. Thus, by \emph{Lemma \ref{Lemma_domain_invariance}}, the mapping ${\bf{\hat \Psi }}_{\upharpoonright \cal U} $ is not injective. As a result, by \emph{Lemma \ref{Lemma_iff_identifiability_condition}}, the vector parameter space ${\boldsymbol{\Theta}}$ is not identifiable.

What's more, according to \emph{Definition \ref{Define_Identifiable_Parameter_Point} and Definition \ref{Define_Identifiable_Parameter_Space}},  the nonidentifiability of the vector parameter space implies that we can find two distinct nonidentifiable points ${\boldsymbol{\theta }}_1 \in {\cal U} \subset {\boldsymbol{\Theta}}$ and ${\boldsymbol{\theta }}_2 \in {\cal U} \subset {\boldsymbol{\Theta}}$ which are observationally equivalent to each other. Note that the set ${\cal U}^*  \buildrel \Delta \over = {\cal U} \backslash \{{\boldsymbol{\theta }}_1, {\boldsymbol{\theta }}_2\} \subset {\cal U}$ is also an open subset of $\boldsymbol{\Theta}$ in ${\mathbbm{R}}^{D_{\boldsymbol{\theta}}}$. Therefore, by the same argument, the restriction ${\bf{\hat \Psi }}_{\upharpoonright {\mathcal U}^*} $ of the mapping ${\bf{\hat \Psi }}$ to ${\mathcal U}^*$ is also not injective, and hence, there also exist two distinct points ${\boldsymbol{\theta }}_1^* \in {\cal U}^* \subset {\cal U} \subset {\boldsymbol{\Theta}}$ and ${\boldsymbol{\theta }}_2^* \in {\cal U}^* \subset {\cal U} \subset {\boldsymbol{\Theta}}$ which are not identifiable. Thus, by induction, there are infinitely many vector parameter points in ${\cal U} \subset {\boldsymbol{\Theta}}$ which are not identifiable. This completes the proof.

\section{Proof of Proposition \ref{Proposition_everypoint_nonidentifiable}} \label{proof_Proposition_everypoint_nonidentifiable}

If the observation $x$ is not quantized, suppose there exist two distinct vector parameter points ${\boldsymbol{\theta}}_1  =  [ \alpha_1, \beta_1]^T$ and ${\boldsymbol{\theta}}_2  =  [ \alpha_2, \beta_2]^T$ which give rise to $f\left( {x\left| {\boldsymbol{\theta}}_1   \right.} \right) = f\left( {x\left| {\boldsymbol{\theta}}_2  \right.} \right)$ for all $x \in {\mathbbm{R}}$, then we can obtain 
\begin{align} \notag
& \ln \frac{{f(x|{{\boldsymbol{\theta }}_1})}}{{f(x|{{\boldsymbol{\theta }}_2})}} \\  \notag
& = \frac{{{\beta _1} - {\beta _2}}}{{2{\beta _1}{\beta _2}}}{x^2} + \left( {\frac{{{\alpha _1}}}{{{\beta _1}}} - \frac{{{\alpha _2}}}{{{\beta _2}}}} \right)x + \frac{1}{2}\left( {\frac{{\alpha _2^2}}{{{\beta _2}}} - \frac{{\alpha _1^2}}{{{\beta _1}}} + \ln \frac{{{\beta _2}}}{{{\beta _1}}}} \right) \\ \label{ln_likelihood_ratio}
&= 0, \; \forall x \in {\mathbbm{R}}.
\end{align}
The fundamental theorem of algebra demonstrates that (\ref{ln_likelihood_ratio}) holds if and only if 
\begin{equation}
\left\{ \begin{array}{l}
\frac{{{\beta _1} - {\beta _2}}}{{2{\beta _1}{\beta _2}}} = 0\\
\frac{{{\alpha _1}}}{{{\beta _1}}} - \frac{{{\alpha _2}}}{{{\beta _2}}} = 0\\
\frac{{\alpha _2^2}}{{{\beta _2}}} - \frac{{\alpha _1^2}}{{{\beta _1}}} + \ln \frac{{{\beta _2}}}{{{\beta _1}}} = 0
\end{array} \right.,
\end{equation}
which implies ${\boldsymbol{\theta}}_1  =  [ \alpha_1, \beta_1^2]^T =    [ \alpha_2, \beta_2^2]^T ={\boldsymbol{\theta}}_2$, and hence, we reach a contradiction. Thus, without quantization, every vector parameter point ${\boldsymbol{\theta}} \in {\boldsymbol{\Theta}}$ is identifiable.

Now, consider the case where the binary quantizer in (\ref{everypoint_nonidentifiable_quantizer}) is employed at the sensor.

For any given vector parameter point ${\boldsymbol{\theta}}_0 \buildrel \Delta \over = {[ {{\alpha _0},\beta _0} ]^T} \in {\boldsymbol{\Theta}}$, let ${{\boldsymbol{\theta }}_\rho } \buildrel \Delta \over = {[ {{\alpha _\rho },\rho \beta _0} ]^T}$ denote a vector parameter point in ${\boldsymbol{\Theta}}$ for some $\rho  \in \left( {0,1} \right)$ and some ${\alpha _\rho }$. We will show that for any $\rho  \in \left( {0,1} \right)$, there exists an ${\alpha _\rho }$ such that ${{\boldsymbol{\theta }}_\rho }$ is observationally equivalent to ${\boldsymbol{\theta}}_0$.

Define a function $g\left( {\alpha ,{\beta }} \right)$ as
\begin{equation} \label{g_function}
g\left( {\alpha ,{\beta }} \right) \buildrel \Delta \over = \Pr \left( {u = 1\left| {\boldsymbol{\theta }} \right.} \right) = \int_a^b {\frac{1}{{\sqrt {2\pi {\beta }} }}{e^{ - \frac{{{{(x - \alpha )}^2}}}{{2{\beta }}}}} dx}. 
\end{equation}
Since $\Pr \left( {u = 2\left| {\boldsymbol{\theta }} \right.} \right) = 1- \Pr \left( {u = 1\left| {\boldsymbol{\theta }} \right.} \right) =  1- g\left( {\alpha ,{\beta }} \right)$, it is clear that if $g\left( {{\alpha _\rho } ,\rho \beta _0} \right) = g\left( {\alpha_0 ,{\beta_0}} \right)$, then by \emph{Definition \ref{Define_Observationally_Equivalent}}, ${{\boldsymbol{\theta }}_\rho }$ are observationally equivalent to ${\boldsymbol{\theta}}_0$, and hence, ${\boldsymbol{\theta }}_0$ is not identifiable.

Since  $I^{(1)}$ and  $I^{(2)}$ are both nonempty sets,
$a$ and $b$ cannot be both unbounded. Without loss of generality, we assume $-\infty \le a < b< \infty$. The case where $ -\infty < a < b \le \infty$ can be proved in a similar way. By (\ref{g_function}) and noting $b < \infty$, we can obtain that for any given ${\alpha _0}$ and $\beta _0$,
\begin{align}  \notag
\mathop {\lim }\limits_{\alpha  \to \infty } g\left( {\alpha ,\rho \beta _0} \right) & = \mathop {\lim }\limits_{\alpha  \to \infty } \int_{\frac{{a - \alpha }}{{\sqrt \rho  \beta }}}^{\frac{{b - \alpha }}{{\sqrt \rho  \beta }}} {\frac{1}{{\sqrt {2\pi } }}{e^{ - \frac{{{x^2}}}{2}}} dx}  \\ \label{g_lower}
& = 0 < g\left( {{\alpha _0},\beta _0} \right).
\end{align}

In the following, we will consider the case where $a = -\infty$ and the case where $a > -\infty$ respectively. We will show that for both cases, there exists some ${{\boldsymbol{\theta }}_\rho } = {[ {{\alpha _\rho },\rho \beta _0} ]^T}$ such that $g\left( {{\alpha _\rho },\rho \beta _0} \right) = g\left( {{\alpha _0},\beta _0} \right)$.

Suppose $a = -\infty$, then for any given ${\alpha _0}$ and $\beta _0$,
\begin{align}  \notag
\mathop {\lim }\limits_{\alpha  \to -\infty } g\left( {\alpha ,\rho \beta _0} \right) & =  \mathop {\lim }\limits_{\alpha  \to  - \infty } \int_{ - \infty }^{\frac{{b - \alpha }}{{\sqrt \rho  \beta }}} {\frac{1}{{\sqrt {2\pi } }}{e^{ - \frac{{{x^2}}}{2}}} dx} \\ \label{g_uper_infty}
&  = 1 > g\left( {{\alpha _0},\beta _0} \right),
\end{align}
Therefore, by (\ref{g_lower}) and (\ref{g_uper_infty}), and noticing that $g\left( {{\alpha  },\beta} \right)$ is a continuous function for all ${\boldsymbol{\theta}}  = {[ {{\alpha },\beta^2} ]^T} \in {\boldsymbol{\Theta}}$, we know that there exists an $\alpha_\rho \in (-\infty, \infty)$ such that
\begin{equation} \label{g_equal_infinite}
g\left( {{\alpha _\rho },\rho \beta _0} \right) = g\left( {{\alpha _0},\beta _0} \right)
\end{equation}
for any given $\rho$ by employing \emph{Intermediate Value Theorem}. 

Suppose $a > -\infty$. Noticing that for any given $\beta$, the equation
\begin{equation}
\frac{\partial }{{\partial \alpha }}g\left( {\alpha ,{\beta }} \right) = \frac{1}{{\sqrt {2\pi {\beta }} }}\left[ {{e^{ - \frac{{{{(a - \alpha )}^2}}}{{2{\beta }}}}} - {e^{ - \frac{{{{(b - \alpha )}^2}}}{{2{\beta }}}}}} \right] = 0
\end{equation}
only admits one solution
\begin{equation}
\alpha  = \frac{1}{2}\left( {a + b} \right) \in \left( { - \infty ,\infty } \right).
\end{equation}
Moreover, since
\begin{equation}
{\left. {\frac{{{\partial ^2}}}{{\partial {\alpha ^2}}}g\left( {\alpha ,\beta } \right)} \right|_{\alpha  = \frac{1}{2}\left( {a + b} \right)}} =  - \frac{{b - a}}{{\sqrt {2\pi {\beta ^3}} }}{e^{ - \frac{{{{(b - a)}^2}}}{{8\beta }}}} < 0,
\end{equation}
we know that $\alpha  = \frac{1}{2}\left( {a + b} \right)$ maximizes the function $g\left( {\alpha ,{\beta }} \right)$ for any given $\beta$. Hence, 
\begin{equation} \label{temp_gaussian}
g\left( {\frac{1}{2}\left( {a + b} \right),\beta _0} \right) \ge g\left( {{\alpha _0},\beta _0} \right).
\end{equation}
Furthermore, note that
\begin{equation} \label{derivative_g_beta}
\frac{\partial }{{\partial \beta }}g\left( {\frac{1}{2}\left( {a + b} \right),\beta } \right) =  - \frac{ {b - a} }{2{\sqrt {2\pi {\beta ^3}} }}{e^{ - \frac{{{{(b - a)}^2}}}{{8\beta }}}} < 0,
\end{equation}
which yields that $g( {\frac{1}{2}( {a + b} ),\beta } )$ is a strictly decreasing function with respect to $\beta$. As a result, by employing (\ref{temp_gaussian}) and (\ref{derivative_g_beta}), we can obtain
\begin{equation} \label{g_upper}
g\left( {\frac{1}{2}\left( {a + b} \right),\rho \beta _0} \right) \!\! > \! g\left( {\frac{1}{2}\left( {a + b} \right),\beta _0} \right) \!\! \ge \! g\left( {{\alpha _0},\beta _0} \right),
\end{equation}
since $\rho \in (0,1)$.
Thus, by (\ref{g_lower}) and (\ref{g_upper}), and by employing \emph{Intermediate Value Theorem}, we know that there exists an $\alpha_\rho \in (\frac{1}{2}( {a + b} ), \infty)$ such that
\begin{equation} \label{g_equal_finite}
g\left( {{\alpha _\rho },\rho \beta _0} \right) = g\left( {{\alpha _0},\beta _0} \right),
\end{equation}
since $g\left( {{\alpha  },\beta} \right)$ is continuous.

By (\ref{g_equal_infinite}) and (\ref{g_equal_finite}), we know that no matter what $a$ and $b$ are, for any given ${\boldsymbol{\theta}}_0  = {[ {{\alpha _0},\beta _0} ]^T} \in {\boldsymbol{\Theta}}$ and for any $\rho \in (0,1)$, there exists some ${{\boldsymbol{\theta }}_\rho } = {[ {{\alpha _\rho },\rho \beta _0} ]^T}$ such that $g\left( {{\alpha _\rho },\rho \beta _0} \right) = g\left( {{\alpha _0},\beta _0} \right)$. Hence, every vector parameter point in $ {\boldsymbol{\Theta}}$ is not identifiable. Moreover, since the set $(0,1)$ is an uncountable set, for any vector parameter point ${\boldsymbol{\theta}}_0 \in {\boldsymbol{\Theta}}$, the set of vector parameter points which are observationally equivalent to ${\boldsymbol{\theta}}_0$ is uncountable. This completes the proof.


\section{Proof of Proposition \ref{Proposition_somepoint_identifiable}} \label{proof_Proposition_somepoint_identifiable}


Define a function $g\left( {\boldsymbol{\theta }} \right)$ as
\begin{align} \notag
g\left( {\boldsymbol{\theta }} \right) & \buildrel \Delta \over = \ln \Pr \left( {u = 1\left| {\boldsymbol{\theta }} \right.} \right) \\ \notag 
& = \ln \int_{{a_1}}^{{b_1}} {\frac{1}{{\sqrt {2\pi } }}{e^{ - \frac{{{{\left( {{x_1} - {\theta _1}} \right)}^2}}}{2}}}d{x_1}}   \\ \label{g_function_vector}
& \qquad 	+ \ln \int_{{a_2}}^{{b_2}} {\frac{1}{{\sqrt {2\pi } }}{e^{ - \frac{{{{\left( {{x_2} - {\theta _2}} \right)}^2}}}{2}}}d{x_2}}.
\end{align}

Note that $(a_1, b_1)$ and $(a_2, b_2)$ are convex sets, and $\frac{1}{{\sqrt {2\pi } }}{e^{ - \frac{{{x^2}}}{2}}}$ is a log-concave function. Hence, $g\left( {\boldsymbol{\theta }} \right)$ is concave, since the integral of a log-concave function over a convex region is log-concave \cite{boyd2004convex}. 

By employing (\ref{g_function_vector}), we can obtain  
\begin{equation}
\frac{d}{{d{\boldsymbol{\theta }}}}g\left( {\boldsymbol{\theta }} \right) =  \left[ \begin{array}{l}
\frac{{{e^{ - \frac{1}{2}{{\left( {{b_1} - {\theta _1}} \right)}^2}}} - {e^{ - \frac{1}{2}{{\left( {{a_1} - {\theta _1}} \right)}^2}}}}}{{\int_{{a_1}}^{{b_1}} {\frac{1}{{2\pi }}{e^{ - \frac{{{{\left( {{x_1} - {\theta _1}} \right)}^2}}}{2}}}d{x_1}} }}\\
\\
\frac{{{e^{ - \frac{1}{2}{{\left( {{b_2} - {\theta _2}} \right)}^2}}} - {e^{ - \frac{1}{2}{{\left( {{a_2} - {\theta _2}} \right)}^2}}}}}{{\int_{{a_2}}^{{b_2}} {\frac{1}{{2\pi }}{e^{ - \frac{{{{\left( {{x_2} - {\theta _2}} \right)}^2}}}{2}}}d{x_2}} }}
\end{array} \right],
\end{equation}
and moreover, by setting $\frac{d}{{d{\boldsymbol{\theta }}}}g\left( {\boldsymbol{\theta }} \right) = {\bf 0}$, we obtain only one solution
\begin{equation}
{{\boldsymbol{\theta }}^ * } = {\left[ {\frac{{{b_1} - {a_1}}}{2},\frac{{{b_2} - {a_2}}}{2}} \right]^T}.
\end{equation}
Thus, $g\left( {\boldsymbol{\theta }} \right)$ achieves the unique globally maximum at ${{\boldsymbol{\theta }}^ * }$, 
since $g\left( {\boldsymbol{\theta }} \right)$ is concave.
Furthermore, since $\Pr \left( {u = 2\left| {\boldsymbol{\theta }} \right.} \right) = 1- \Pr \left( {u = 1\left| {\boldsymbol{\theta }} \right.} \right)$, it is clear that if there exists a vector parameter point ${\boldsymbol{\theta}}$ such that $g\left( {\boldsymbol{\theta }} \right) \ne g\left( {\boldsymbol{\theta }}' \right) $ for all ${\boldsymbol{\theta }}' \in {\boldsymbol{\Theta}}\backslash \{{\boldsymbol{\theta }}\}$, then by \emph{Definition \ref{Define_Identifiable_Parameter_Point}}, ${\boldsymbol{\theta}}$ is identifiable. As a result, ${\boldsymbol{\theta }}^* $ is an identifiable vector parameter point in 
$\boldsymbol{\Theta}$. This completes the proof.

%
%

\bibliographystyle{IEEEtran}
\bibliography{Attack}

\end{document}